\documentclass[12pt]{article}
\usepackage{fancybox,%
            fancyheadings,%
            epsf,%
            theorem,%
            eepic,%
            psfig,%
            epsfig,%
            amssymb,%
            amsmath,%
            pifont}
\usepackage%
[
nodraft,
lettersize
]
{symbols}
 \norems                    
\newcommand{\wt}{\widetilde}
\renewcommand{\putfigtexh}[4][\figpath]{
\let\auxfigpath=\figpath \let\figpath=#1
\begin{figure}[h]            
\let\auxgraphics=\includegraphics
\renewcommand{\includegraphics}[1]{\auxgraphics{\figpath/#2.pstex}}
 \vbox{ 
      \centerline{ \input{\figpath/#2.pstex_t}}
      }
\caption{#4}
\label{#3}                   
\let\includegraphics=\auxgraphics
\end{figure}                    
\let\figpath=\auxfigpath 
}

 \let\c=\cedille
                  
 \def\coauthoronesname{D}
 
 \def\coauthortwoname{A}

\renewcommand\norm[1]{\left|#1\right|}
\newcommand\Norm[1]{\left\|#1\right\|}
\def\HdD{{\cal H}(\delta,\Delta)}

\def\BD{{\cal B}(\Delta)}

\def\Bmu{{\cal B}(\mu)}
\def\wtg#1#2{{\widetilde g^{\ \! #1}_{#2}}}

\begin{document}

\title{Matrosov's theorem using a family of auxiliary functions: \\
an analysis tool to aid time-varying nonlinear control design\footnote{This work is supported in part by a CNRS-NSF collaboration project and was partially done while the first two authors were visiting the CCEC at UCSB. It is also supported by the AFOSR
under grant F49620-00-1-0106 and the NSF under grants ECS-9988813 and INT-9910030. }}
\author{Andrew R. Teel$^\bullet$ \quad %
Antonio Lor\'{\i}a$^\dagger$\quad %
Elena Panteley$^*$ \quad %
Dobrivoje Popovi\'c$^\bullet$ \quad %
\address{%
$^\bullet$Center for Control Engineering and Computation (CCEC),  Department of Electrical and Computer Engineering, University  of California,  Santa  Barbara, CA 93106-9560, USA
\\[2mm]
$^\dagger$C.N.R.S, UMR 8506, Laboratoire de Signaux et Syst\`emes,
Sup\'elec, Plateau de Moulon, 91192 Gif s/Yvette, France
\\[2mm]
$^*$ IPME, Acad. of Sciences of Russia
 61, Bolshoi Ave,
 St. Petersburg, Russia
 elena@ieee.org
}}

\date{ }

\sloppy
\maketitle

\begin{abstract} 
We present a new result on uniform attractivity of the origin for
nonlinear time-varying systems.  Our theorem generalizes Matrosov's
theorem which extends, in a certain manner, Krasovsk\u{\i}-LaSalle
invariance principle to the case of general nonlinear time-varying
systems. We show the utility of our theorem by addressing a control problem of port interconnected driftless systems. The latter includes as special case, the control of chained-form nonholonomic systems which has been extensively studied 
in the literature.
\end{abstract}


\section{Introduction}

Many nonlinear control algorithms rely heavily on analysis tools that
establish convergence to the origin for trajectories of a time-varying
nonlinear system having a uniformly stable origin.  For time-invariant
problems, the typical analysis tool used is the Krasovsk\u{\i}-LaSalle
invariance principle.  It is the key result that leads to the
so-called Jurdjevi\'c-Quinn control algorithm for open-loop stable
nonlinear control systems.  When the closed-loop is time-varying, one
tool that is often used is Barb\u{a}lat's Lemma \cite{BARBALAT}.  In adaptive control,
Barb\u{a}lat's Lemma is frequently relied upon to establish convergence to
zero of part of the state.  Barb\u{a}lat's Lemma has also been used to
establish convergence to the origin for a class of nonholonomic
systems controlled by smooth time-varying feedback.  For time-varying
systems, another tool that has been used, but more sparingly, is
Matrosov's theorem \cite{MAT62,ROUMAW}.  It was used e.g. in \cite{PADPAN} to establish one of the first results on uniform global asymptotic stability (UGAS) of
robot manipulators in closed loop with a tracking controller. It also
appears in the context of adaptive control in \cite{MARTOM} and output feedback control in \cite{NICTOM1}.

Matrosov's theorem first appeared in \cite{MAT62}. It pertains to the
situation where one has a continuously differentiable ($C^{1}$)
Lyapunov function that establishes uniform stability and also a
$C^{1}$ auxiliary function with appropriate properties.  In
particular, the auxiliary function should be bounded uniformly in time
on bounded regions of the state space, and should have a ``definitely
nonzero'' derivative on the set where a given, continuous,
time-independent non-positive upper bound on the Lyapunov function's
derivative vanishes.  To illustrate how Matrosov's theorem works on tracking control problems it is worth looking at the example of tracking control of robots. 

\begin{example}\label{ex:intro}
Let us consider the Lagrangian model of a rigid-joints robot  (see e.g. \cite{SPOVID})
\[
D(q)\ddot q + C(q,\dot q)\dot q + g(q) = u
\]
where the inertia matrix $D$ is positive definite for all $q\in\mR^n$, $\dot D - 2C$ is skew-symmetric and $u\in\mR^n$ are control torques. The control problem is to make the robot follow a smooth reference trajectory  $\qd(t)$ such that $\max\{\norm{\qd},\,\norm{\dqd},\,\norm{\ddqd}\}\leq \beta_d$. For this, we apply the control law (cf. \cite{PADPAN})
\[
u = D(q)\ddqd + C(q,\dot q)\dqd + g(q) - k_p\tq - k_d\dtq
\]
where $\tq:=q-\qd$, $k_p$, $k_d>0$. Now we would like to analyze the stability of the closed loop system, so we use the energy-based Lyapunov function 
\[
V= \frac{1}{2} \left(\dtq^\top D(q)\dtq +  k_p \norm{\tq}^2 \right)
\]
to obtain that $\dot V = - k_d \norm{\dtq}^2 \leq 0$ and hence, that the origin of the system is uniformly globally stable (UGS). From this equality, it is also fairly standard to invoke Barb\u{a}lat's Lemma to conclude that $\dtq\to 0$ (see e.g. \cite{ORTSPO}). Alternatively, we can use the auxiliary function $W = \tq^\top D \dtq$ which is uniformly bounded on compact sets of the state and whose total time derivative satisfies\footnote{Using the facts that $\norm{D(q)}\leq d_M$ and $\norm{C(q,\dq)} \leq k_c\norm{\dq}$. See \cite{SPOVID} for further details.} 
$$
\dot W \leq  - k_p\norm{\tq}^2 + \norm{\dtq}[\,d_M \norm{\dtq} + k_d\norm{\tq} + k_c(\beta_d + \norm{\dtq})\norm{\tq}\,] =: Y\,.
$$
Intuitively, since we know that $\dtq \to 0$ we may think that $Y \to - k_p\norm{ \tq}^2$ so roughly, if we wait long enough (and we can do that because the system is UGS) we will have that for large $t$,  $\dot W(t) \leq  - k_p\norm{ \tq(t)}^2$. In particular, $\tq$ is square-integrable so invoking again Barb\u{a}lat's Lemma one may accept that it should hold that $\tq\to 0$ as well. 

This idea can be made rigorous to actually prove UGAS via Matrosov's theorem: one needs to observe that $Y$ is sign-definite (it is actually negative for all nonzero values of the state) on the set $\{\dot V = 0\}$ that is, on  $\{\dtq = 0, \, \tq\in\mR^n\}$. The argument at work here is roughly speaking that the sign-definiteness of $Y$ and the fact that $W=0$ on the set $\{\dot V = 0\}$ imply that the system's trajectories cannot remain trapped on the set $\{\dtq = 0, \, \tq\in\mR^n\}$ unless they go to zero.  See \cite[p. 263]{HAH} for a rigorous  development on this idea and \cite{PADPAN} for a rigorous analysis of this control system based on Matrosov's theorem. 
\end{example}

Thus, Matrosov's theorem can be regarded, to some extent, as an invariance principle for non-autonomous systems. There have been several interesting extensions of Matrosov's theorem over the years, mostly found in the work of \cite{ROUMAW,ROUHABLAL} and references therein. In \cite{ROU-MAT}, the author considers a vector of auxiliary functions, while in \cite{CORROU} the authors consider a family of (possibly uncountable) auxiliary functions.  In both cases, the behavior of the auxiliary functions are referenced to the set where the upper bound on the derivative of the Lyapunov function vanishes. In \cite{CORROU}, extensions are given that pertain to stability of sets and that allow locally Lipschitz Lyapunov and auxiliary functions.  In \cite{PADPAN} a simplified condition is given for checking that the auxiliary function has a sign-definite derivative.  

Most recently, in \cite{INTLEMMCSS} the authors have extended Matrosov's
theorem to pertain to differential inclusions, at the same time
addressing stability of sets, using locally Lipschitz auxiliary
functions and weakening the requirements on the upper bound of the
derivative of the Lyapunov function.

\subsection{Contributions of this paper}

In this paper we  develop a version of Matrosov's theorem that uses a finite family of auxiliary functions to establish uniform convergence.  For simplicity of exposition, we will limit our discussion to differential equations, uniform asymptotic stability of the origin, and a means of checking a sign-definiteness condition that is similar to what was used in \cite{PADPAN}.  Nevertheless, the results that we present here extend easily to more general settings (stability of sets, differential inclusions, $\varepsilon - \delta$ characterizations of sign-definiteness, etc.)

Since we use a finite family of auxiliary functions, it is natural to wonder about the comparison to the vector of auxiliary functions used in \cite{ROU-MAT} and the uncountable family of auxiliary functions used in \cite{CORROU}.  The main difference is that, whereas in the above references the behavior of all of the auxiliary functions and their derivatives is referenced to the set where the upper bound on the derivative of the Lyapunov function vanishes, our auxiliary functions are ordered and the behavior of an auxiliary function and its derivative is referenced to the set where the upper bounds on the derivatives of all of the preceding auxiliary functions vanish.  We assume uniform stability, rather than assuming that we have a Lyapunov function which establishes it, and then none of our auxiliary functions is assumed to be sign definite.  

Finally, we remark that the motivation to  have a family of auxiliary functions as opposed to 2 functions as in the case of the original Matrosov's theorem, should become clearer through the case-studies in control design that we will address: control of port-controlled driftless systems, and as a particular case, chained-form and ``skew-symmetric'' systems (see e.g. \cite{MURSAS,SAM3} respectively). To anticipate the development of these control applications we shall say that a non-obvious fundamental property of the mechanical system in Example \ref{ex:intro} that makes Matrosov's theorem work to prove UGAS is that the system has relative degree one from the {\em passive} output $\dtq$. This property still holds for chained form systems of $3$ states however, it is lost for higher-dimensional systems. While there seems to be no precise arithmetic correlation between the degree of the system and the number of functions needed, it shall become apparent from the case-studies, that it is natural to have a number of auxiliary functions which is directly related to the number of states in the chain of (nonlinear) integrators and therefore, to the relative degree of the system.

Thus, we believe that the analysis tool that we are presenting may become an efficient tool to aid time-varying nonlinear control design, and may contribute to making the idea of using auxiliary functions for nonlinear systems, as introduced by Matrosov, a more versatile concept.

The rest of the paper is organized as follows. For clarity of exposition, we present our results from the most particular to the most general: we start in Section \ref{sec:nonhol3states} with the ``simple'' control problem of stabilizing a 3-states chained-form system. The way we solve the problem of analyzing the stability of this system in closed loop with previously reported control laws, will motivate the sufficient conditions of our main Theorem. The latter is presented in Section \ref{sec:matrosov}. In Section \ref{sec:channels} we present our main result in control design which covers chained-form systems as a particular case.  The proof of our main theorem is presented in Section \ref{sec:proof} and we conclude with some remarks in Section \ref{sec:concl}.

\subsection{Mathematical preliminaries} 
We use $\norm{\cdot}$ to denote the Euclidean norm of vectors and the induced norm of matrices. We use $\Norm{\cdot}_p$, where $p\in [1,\infty]$, denotes the ${\cal L}_p$ norm of time signals.  In particular,  for $t_\circ\in\mR$ and a measurable function $\phi: \mR_{ \geq t_{\circ} } \rightarrow \mR^{n}$, by $\Norm{\phi}_p$ we mean $(\int_{\t\circ}^\infty \norm{\phi(t)}^p dt )^{1/p}$ for $p\in [1,\infty)$ and $\Norm{\phi}_\infty$ denotes the quantity $\mbox{\rm ess sup}_{t\geq \t\circ} \norm{\phi(t)}$. For two constants $\Delta\geq\delta\geq 0$ we define ${\cal H}(\delta,\Delta) \defeq \{x\in \mR^n \, : \, \delta \leq \norm{x} \leq \Delta\}$.  We also will use ${\cal B}(r):={\cal H}(0,r)$. A continuous function $\rho:\mRp\to \mRp$ is of class ${\cal N}$ if it is non-decreasing. A function $\gamma : \mRp \to \mRp$ is of class ${\cal K}$ ($\gamma \in {\cal K}$), if it is continuous, strictly increasing and zero at zero; $\gamma \in {\cal K}_\infty$ if in addition, $\gamma(s) \to \infty$ as $s\to \infty$.  A function $\alpha: \mRp \to \mRp$ is of class ${\mathcal P}{\mathcal D}$ if it is continuous and positive definite. We denote by $x(\cdot,\t\circ,\x\circ)$, the solutions of the differential equation
\begin{equation}
\label{eq:nltv}
 \dot x = f(t,x)
\end{equation}
with initial conditions $(\t\circ,\x\circ)$. For a function $V:\mR\times\mR^n\to \mR$ we define its derivative in the direction $(1,f(t,x)^{T})^{T}$ as $\dot V(t,x):= \jac{V}{t}+ \jac{V}{x}f(t,x)$\,. This abuse of notation is reasonable because of (\ref{eq:nltv}). Furthermore, with an abuse of notation we will use the same definition to express the derivative of locally Lipschitz functions. For the latter, $\dot V(t,x)$ is defined everywhere except on a set of measure zero where the gradient of $V$ is not defined, i.e., almost everywhere. When clear from the context and to simplify the notation we will also use $\nabla f(x)$ to define $\jac{f(x)}{x}$.

\begin{defin}[Uniform global stability]
  \label{def:stability}
  The origin of the system \rref{eq:nltv} is said to be uniformly 
  globally stable (UGS) if there exists $\gamma\in {\cal K}_\infty$ such that, for
  each $(t_{\circ},\x\circ) \in \mR \times \mR^{n}$ each solution $x(\cdot,t_{\circ},x_{\circ})$ satisfies 
\begin{equation}
\label{us}
  \norm{x(t,t_{\circ},x_{\circ})} \leq \gamma(\norm{x_{\circ}})  \qquad \forall \; t \geq t_{\circ} \ .
\end{equation}
\end{defin}
\begin{defin}[Uniform global attractivity]
The origin of the system (\ref{eq:nltv}) is said to be  uniformly globally attractive  (UGA) if  for each  $r,\,\sigma>0$  there exists $T>0$ such that
\begin{equation}
\label{UA}
\norm{x_{\circ}}\leq r  \, \Longrightarrow \, \|x(t,t_{\circ},x_{\circ})\|
\leq \sigma \qquad  \forall\, t \geq t_{\circ} + T\,
  \ .
\end{equation}
\end{defin}
Furthermore, we say that the (origin of the) system  is uniformly globally asymptotically stable (UGAS) if it is UGS and UGA.

The results on control design that we will present rely on the notion of uniform $\delta$-persistency of excitation (U$\delta$-PE) which is defined next for completeness and ease of reference.
\begin{definition}[see \cite{DPECDC02}]\label{def:uped}
A function  $\phi(\cdot,\cdot)$ where $t\mapsto\phi(t,x)$ is locally integrable, is said to be uniformly $\delta$-persistently exciting (U\ped) with respect to $x_1$ if for each $x \in \mR^n\backslash \{x_1=0\} $ there exist $\delta >0$, $T>0$ and $\mu>0$ 
s.t. $\forall t \in \mR$,
\begin{equation}\label{peofpsi}
  \norm{z-x} \leq \delta  \quad \Longrightarrow \quad \int_{t}^{t+T} \norm{\phi(\tau,z)} d\tau \geq \mu\,.
\end{equation}
\end{definition}
Roughly speaking this property means that the function $u(t,x(t))$ is persistently exciting\footnote{We recall that  in the  adaptive control literature a function $t\mapsto \psi(t)$ is called PE if there exist $\mu$, $T>0$ such that 
\[
\int_t^{t+T} \psi(\tau)^2d\tau \geq \mu \quad \forall\, t\geq 0\,.
 \]
} (PE) whenever the trajectories $x_1(t)$ are away from a $\delta$-neighborhood of the origin. We will also make use of the following observations.
\begin{fact}\label{fact:prodphis}
For locally Lipschitz functions, uniformly in $t$, and such that $\phi_i(t,0)\equiv 0$ we have that if the product $\prod_{i=1}^n\phi_i(t,x)$ is U\ped\ then, necessarily each function $\phi_i(t,x)$  is U\ped. In particular, $\phi_i(t,x)^p$ is U\ped\ for any integer $p$. 
\end{fact}
\begin{property}\label{prop:filter}
Consider the differential equation
\begin{equation}
\label{eq:dphif}
\dot \Phi_f = -a \Phi_f + \phi(t,z),\qquad a >0
\end{equation}
where  $\phi(t,x)$ is U\ped\ with respect to $z_1$. Then, defining $\tilde \phi_f(\cdot\,,t_{f_\circ},\phi_{f_\circ},x)$ as the solution of \rref{eq:dphif}, $x_2:=\col[t_{f_\circ},\phi_{f_\circ},z_2]$ and $x_1:=z_1$, the function $\phi_f(t,x):=\tilde\phi_f(t,t_{f_\circ},\phi_{f_\circ},z)$ is U\ped\ with respect to $x_1$ 
\end{property}
\begin{fact}
\label{fact:onpeofpsi}
We have from \cite[Lemma 3]{DPECDC02} that there exist a continuous function $\theta_\Delta: \mR_{>0}\to\mR_{>0}$ and $\gamma_\Delta\in\cK$ such that the lower bound on the integral in \rref{peofpsi} holds for all $x_2\in{\cal B}(\Delta)\backslash \{x_2=0\}$ with $T=\theta_\Delta(\norm{x_2})$ and $\mu=\gamma_\Delta(\norm{x_2})$.
\end{fact}

\section{The nonholonomic integrator: a case study}
\label{sec:nonhol3states}

 In order to motivate our extended Matrosov's theorem (cf. Theorem \ref{thm:main:new:mat}) and to illustrate the formulation of the sufficient conditions that it imposes for UGAS, we find it appropriate to deal with a well known control problem: the stabilization of a 3-states chain form system (see e.g. \cite{MURSAS}) by smooth time-varying feedback. 

The problem of stabilization of nonholonomic systems of any dimension is well documented and numerous results have been proposed in the literature. We do not intend to revisit this problem in detail and therefore, we suggest that interested readers see for instance \cite{KOLMCC} for a tutorial with a very complete literature review. See also the more recent references \cite{SAMSIAM02,ERJENSTHESIS}.

Consider the problem of stabilizing the nonholonomic {\em chained} system 
\begin{subequations}\label{eq:nonhol}
\begin{eqnarray}
\dot x_1 &=& u(t,x)\\
\dot x_2 &=& u(t,x)x_3\\
\dot x_3 &=& v(t,x)
\end{eqnarray}
\end{subequations}
by  smooth feedback controls $u$ and $v$. In particular, we will revisit the controllers proposed in  \cite{SAM3,uped} and present a new proof\footnote{These controllers were originally proposed in \cite{SAM3} but uniformity was proved only in \cite{uped}.} of UGAS of the closed loop. The result can be proved for any number of states as we will see later but for the motivation of the conditions in Theorem \ref{thm:main:new:mat}, three states are enough. 

 The approach presented below  encompasses in particular  periodic time-varying feedbacks as considered for instance in  \cite{JIANONHOL96}. With regard to this reference it is also interesting to observe that the author used Krasovsk\u{\i}-LaSalle invariance principle to obtain a {\em direct} proof of global asymptotic stability. Here, we use our generalization of Matrosov's theorem which may be considered as the extension of the invariance principle, for general nonautonomous systems.   

We start by recalling the control laws from \cite{SAM3,uped} which are 
\begin{subequations}\label{control}
\begin{eqnarray}
\label{control:a} v(t,x) & =& - x_3 - u(t,x)x_2\\
\label{control:b} u(t,x) &= &-x_1 + h(t,x_{2,3})\,, 
\end{eqnarray}
\end{subequations}
where $\ x_{2,3}:= \col[x_2,\, x_3]$, used in the cited references and where $h(\cdot,\cdot)$ is a smooth function defined below. 
\begin{proposition}
\label{prop:nonhol}
Consider the system \rref{eq:nonhol} in closed loop with \rref{control}. Let
\begin{assumption*}
\label{ass:smoothness}
the map $(t,x_{2,3})\mapsto h(t,x_{2,3})$ be such that $h(t,0)\equiv 0$, all its first and second partial derivatives are uniformly bounded by $\rho(\norm{x_{2,3}})$ where $\rho(\cdot)$ is a non decreasing function\footnote{For simplicity, we will $\rho(\cdot)$ for a generic bound on any function which is uniformly bounded in $t$.}
\end{assumption*}
and, defining
\begin{equation}
\label{defofpsi}
\psi(t,x_2):= \jac{h_\circ}{t}\,, \qquad 
h_\circ(t,x_2) := h\left(t,\begin{bmatrix} x_2\\ 0 \end{bmatrix}\right)
\end{equation}
assume that $\psi(t,x_2)$ is U\ped\ with respect to $x_2$. Then, the origin of the closed loop system is UGAS.
\end{proposition}
The proof of this claim under a similar property to that defined by \rref{peofpsi} was first presented in \cite{uped} relying on a ``cascades'' argument and linear time varying systems theory. The proof that we propose here follows directly  by constructing certain functions which allow to show that all signals are ${\cal L}_p$-integrable. This will motivate the conditions of our main theorem. Moreover, this proof can be extended for more than three states. In contrast to this, the only  proof of UGAS with similar controls that we know for systems with $n>3$ states requires the application of numerous auxiliary results for linear time-varying systems (see \cite{NONHOLEJC} and Section \ref{sec:exnonhol}). 

The general intuition to establish the proof can be explained in terms of Krasovsk\u{\i}-LaSalle invariance principle. For this, let us restrict our attention to periodic feedbacks (as in \cite{JIANONHOL96}). The closed loop system of \rref{eq:nonhol} with \rref{control} is 
\begin{subequations}\label{cl3states}
\begin{eqnarray}
\label{cl3states:a}
\dot x_1 & = & -x_1 + h(t,x_{2,3})\\ 
\label{cl3states:b}
\dot x_2 & = & u(t,x)x_3 \\
\label{cl3states:c}
\dot x_3 & = & -x_3 + u(t,x)x_2\,.
\end{eqnarray}
\end{subequations}
Firstly, taking the derivative of $V_1(t,x):= \frac{1}{2} \norm{x_{2,3}}^2$ we obtain that
\begin{equation}\label{mm}
 \dot V_1(t,x)\leq -x_3^2 \leq 0\,.
\end{equation}
From this inequality we obtain that $x_3\to 0$ hence, we may also admit from \rref{cl3states:c} that $u(t,x)x_2\to 0$. In addition to this, we have from \rref{cl3states:b} that $x_2\to $ const. This means that $u(t,x)$ tends to a steady-state value which we denote by $\omega(t,x_2)$. Now, since $\omega(t,x_2)$ is periodic in $t$ it is reasonable to assume that  it is also  sufficiently rich (persistently exciting) for each $x_2 \neq 0$. If this is the case, then from the conjecture that $\omega(t,x_2)x_2\to 0$ necessarily the only constant value which $x_2$ may converge to is zero. Finally, the convergence of $x_1$ is obtained from \rref{cl3states:a} and the fact that $h(t,0)\equiv 0$. 

The clear drawback of such an argument is that it cannot be made precise for general non-autonomous systems since it relies on Krasovsk\u{\i}-LaSalle invariance principle. We establish below a rigorous proof based on the intuiton developed so far. 

\subsection{Proof of Proposition \ref{prop:nonhol}}
\label{sec:proof3states}
Firstly, UGS follows from \rref{mm} and the following arguments. Integrating \rref{mm} from $t_\circ$ to $\infty$ we obtain that $\norm{x_{2,3}(t)} \leq \norm{x_{{2,3}_\circ}}$ for all $t\geq t_\circ$. Technically, this inequality is valid only on the interval of existence of the solutions. Integrating on this window and using the fact that $\norm{x_{2,3}(t)}$ is bounded on the maximal interval of definition we proceed to integrate the $x_1$-equation \rref{cl3states:a} to obtain that $\norm{x_1(t)} \leq \norm{x_{1_\circ}} + \rho(\norm{x_{{2,3}_\circ}})$ for all $t\geq t_\circ$ and where $\rho(\norm{x_{{2,3}_\circ}})$ is a generic bound on $h(\cdot,x_{2,3})$ which exists due to Assumption \ref{ass:smoothness}. Therefore, the solutions exist for all $t$ and actually, the origin is UGS. Notice also that {\em here}, $\rho(0)=0$ since $h(t,0)\equiv 0$.

To prove attractivity we consider other differentiable functions which are bounded and have bounded derivatives on balls 
$
{\cal B}(\Delta) 
$. Our starting point in the pursuit of these additional functions to combine with $V_1$ is the observation that any terms in the derivative of subsequent auxiliary functions that vanish with $x_3$, can be ignored since we know that\footnote{The precise argument that makes this rough observation hold true is contained in the upcoming Theorem \ref{thm:main:new:mat}.} $x_3\to 0$. So, for example, we can take 
\begin{equation}\label{v2}
V_2(t,x):= x_3 u(t,x) x_2
\end{equation}
and, defining $\phi(t,x):= u(t,x)x_2$, we obtain that 
\[
\dot V_2(t,x) \leq -\phi(t,x)^2 - x_3\phi(t,x) + \dot u(t,x)x_2x_3 + [x_3u(t,x)]^2\,.
\]
From smoothness of $u(\cdot,\cdot)$ and Assumption \ref{ass:smoothness} we obtain that its total derivative is bounded for bounded $x$, uniformly in $t$. 
 In the sequel we will use the number $\nu>0$ as a generic bound on continuous functions over compact sets. With this under consideration, we have that
\[
\dot V_2(t,x) \leq -\phi(t,x)^{2} + \nu\norm{x_3}\qquad \forall\,(t,x)\in \mR\times {\cal B}(\Delta)\,. 
\]
For the sequel, we see that we can now  ignore as well  terms in derivatives of auxiliary functions, that vanish with $\phi(t,x)= u(t,x)x_2$. As a matter of fact, from \rref{eq:nonhol} and \rref{control}, we now see that we can ignore also terms that vanish with $x_3$, $\dot x_3$ and $\dot x_2$. What is more, if we were  to consider the dynamics of the closed loop system system when  $x_3 \equiv 0$ and $x_2$ is constant the $x_1$-equation would define the dynamics of a linear system with a time varying input parameterized by a $constant$ $\bar x_2$, i.e,
\[
\dot{\bar x}_1 = -\bar x_1 + h_\circ(t,\bar x_2) \ :=\ \bar u(t,\bar x)\,.
\]
Differentiating on both sides and owing to the fact that $\bar x_2\equiv$const, we obtain that $\bar u$ satisfies the differential equation 
$\dot {\bar u} = -\bar u + \psi(t,\bar x_2)\,$
and whose solution is 
\begin{equation}\label{omegadef}
\omega(t,\bar x_2)  :=  \int_{-\infty}^t \mbox{e}^{-(t-\tau)} \psi(\tau,\bar x_2)d\tau\,.
\end{equation}
Based on these observations, we will  introduce the next function  with the aim at concluding something about the difference between $\dot x_1(\cdot)$ and the steady state solution $\omega(\cdot,\bar x_2)$. That is, we define our third auxiliary function as $V_3(t,x):= \zeta(t,x)^2$, where
\begin{equation}\label{zetadef}
\zeta(t,x)  := x_1 - h_\circ(t,x_2) + \omega(t,x_2) \,,
\end{equation}
and observe that the time derivative of $V_3(t,x)$ along the trajectories of\footnote{We emphasize that the definition of the function $V_3(\cdot,\cdot)$ is motivated by the behavior of the system on the manifold $\{ x_3\equiv 0,\,x_2\equiv \mbox{const}\}$ however, we consider its total derivative along the trajectories of the closed loop system, on ${\cal B}(\Delta)$. In other words, we {\em do not} analyze the system's dynamics {\em only} on the defined manifold. }
\begin{eqnarray*}
\dot \zeta(t,x) & \!\!\! = & \!\!\!- \zeta(t,x) + h(t,x_{2,3}) - h_\circ(t,x_2) - \jac{h_\circ}{x_2} u(t,x)x_3  
+ \left(\int_{-\infty}^t \mbox{e}^{-(t-\tau)}\jac{\psi}{x_2}\right)u(t,x)x_3\,
\end{eqnarray*}
yields 
\begin{equation}
\label{dotv3}
\dot V_3(t,x) \leq -\zeta(t,x)^2 + \nu\norm{x_3}\,.
\end{equation}
To obtain this inequality  we have  used the smoothness of all functions, Assumption \ref{ass:smoothness}  (which in particular implies that  $h(t,\cdot)$ is locally Lipschitz, uniformly in $t$), and the compactness of ${\cal B}(\Delta)$.

Our fourth function is introduced to be used in combination with the U\ped\ property of $\psi(t,x)$ in order to infer that the only constant value that $x_2(\cdot)$ may converge to is zero.  This function is
\begin{equation}
\label{v4}
 V_{4}(t,x) := - \int_{t}^{\infty} \mbox{e}^{(t-\tau)} \norm {\omega(\tau,x_{2})x_2}^2d\tau\,
\end{equation}
which satisfies for any $T>0$,
\begin{equation}
  V_{4}(t,x) \leq - e^{-T} \int_{t}^{t+T} \norm {\omega(\tau,x_{2})x_2}^2d\tau\,
\end{equation}
and we claim that from the above and \rref{peofpsi} it follows that there exists  a continuous, non-decreasing function $\gamma(\cdot)$ such that $\gamma(0)=0$ and for all $x_2\in{\cal B}( \Delta)$, 
\begin{equation}
\label{31}
  V_{4}(t,x) \leq - \gamma(\norm{x_2})\norm{x_2}^2\,.
\end{equation}
This follows by appealing to Property \ref{prop:filter} of U\ped\ functions and observing that $\omega(t,x_2)$ is defined by the differential equation 
$$
\dot \omega = -\omega + \psi(t,x_2)\,.
$$
Since $\psi(t,x_2)$ is U\ped\ with repsect to $x_2$, it follows that there exists  $\mu_*>0$ such that $\omega(t,x_2)$ satisfies a bound like \rref{peofpsi}. Furthermore, appealing to Fact \ref{fact:onpeofpsi}\ and using the inequality $\int_t^{t+T}\norm{f(\tau)}d\tau \leq (T\int_t^{t+T}\norm{f(\tau)}^2d\tau )^{1/2}$ the claim follows with $\gamma(s):= \min\left\{ s,\, \frac{e^{- \theta_\Delta(s)}}{\theta_\Delta(s)}\gamma_\Delta(s)^2\right\}$.

We proceed now to evaluate the time derivative of $V_4(t,x)$ along the trajectories of the closed loop system.  To that end, we write 
\begin{eqnarray}
\label{dv4dx2}
\jac{V_4}{x_2} & \!\!\! = & \!\!\! -\int^{\infty}_t \mbox{e}^{(t-\tau)} 2\omega(\tau,x_2)x_2
\left[\jac{\omega}{x_2}x_2 + \omega(\tau,x_2)\right]d\tau\nonumber \\ && \\
\label{domdx2}
\jac{\omega}{x_2} & \!\!\! = & \!\!\! \int_{-\infty}^t\mbox{e}^{-(t-\tau)}\jac{}{x_2}\left(\jac{h_\circ}{t}\right)d\tau\\
\jac{V_4}{t} & \!\!\! = & \!\!\! \mbox{e}^{(t-\tau)} \norm{\omega(\tau,x_2)x_2}^2 \big|_{\tau=t} 
- \int^{\infty}_t \jac{}{t}\left[ \mbox{e}^{(t-\tau)}\norm{\omega(\tau,x_2)x_2}^2\right] d\tau
\end{eqnarray}
and observe that due to Assumption \ref{ass:smoothness} all the partial derivatives in \rref{dv4dx2} and \rref{domdx2} which, are functions of $(t,x_2)$,  are uniformly bounded in $t$ by a generic bound that we denote $\rho(\norm{x_2})$. Since $u(\cdot,\cdot)$ also satisfies this boundedness property with $\rho(\norm{x})$ we finally obtain using \rref{zetadef}, that for all $t\in\mR$ and all $x\in {\cal B}(\Delta)$, 
\begin{eqnarray*}
\dot V_4(t,x) &\leq & V_4(t,x) + \norm{\omega(t,x_2)x_2}^2 + \nu\norm{x_3}\\
&& = V_4(t,x) + | \,x_2[\zeta(t,x)+u(t,x)  
- h(t,x_{2,3}) + h_\circ(t,x_2)]\, |^2 + \nu\norm{x_3}
\end{eqnarray*}
and using \rref{31} and the fact that $\norm{h_\circ(t,x_2)- h(t,x_{2,3})}\leq \nu\norm{x_3}$ for all $t\in\mR$ and $x\in{\cal B}(\Delta)$ we finally obtain that
\[
\dot V_4(t,x) \leq -\gamma(\norm{x_2}) \norm{x_2}^2 + \nu[\, \zeta(t,x)^2 + \phi(t,x)^2  + \norm{x_3}\, ] \,.
\]
Thus, $V_4(t,x)$ helps us to see that for the subsequent functions we can also ignore all the terms vanishing with $x_2$. It is only left to find a function whose derivative is bounded by a negative term of $\norm{x_1}$ and possibly positive terms of $\norm{x_2}$, $\norm{x_3}$. For this, we introduce $V_5(t,x)= x_1^2$ whose total time derivative yields
\[
\dot V_5(t,x) \leq -x_1^2 + 2\norm{x_1}\norm{h(t,x_{2,3})} \leq -x_1^2 + \nu\norm{x_{2,3}}\,.
\]

Summarizing, we have that the functions 
\begin{eqnarray}
V_1(t,x)&:=& \frac{1}{2}\left(x_2^2 + x_3^2\right)\\
V_2(t,x)&:=&  x_3 u(t,x)x_2\\
V_3(t,x)&:=& \zeta(t,x)^2\qquad \mbox{\small  see \rref{zetadef} and \rref{omegadef} }\\
V_4(t,x)&:=& - \int_{t}^{\infty} \mbox{e}^{(t-\tau)} \norm {\omega(\tau,x_{2})x_2}^2d\tau\,\\
V_5(t,x)&:=&  x_1^2 \,.
\end{eqnarray}
satisfy that $\dot V_i(t,x)\leq Y_i(t,x)$ for almost all $(t,x)\in \mR \times {\cal B}(\Delta)$ with  $\phi_1(t,x):=\zeta(t,x)$, $\phi_2(t,x):=\phi(t,x)$, and
\begin{eqnarray*}
Y_1(t,x)&:=& -x_3^2\\
Y_2(t,x)&:=&  -\phi_2(t,x)^2 + \nu\norm{x_3} \\
Y_3(t,x)&:=& -\phi_1(t,x)^2 + \nu\norm{x_3} \\
Y_4(t,x)&:=& -\gamma(\norm{x_2}) \norm{x_2}^2 + \nu[\, \phi_1(t,x)^2 + \phi_2(t,x)^2  + 
\norm{x_3}\, ] \\
Y_5(t,x)&:=&  -x_1^2 + \nu\norm{x_{2,3}}\,.\\[-10mm]
\end{eqnarray*}

To conclude (non-uniform) convergence we may appeal successively to Barb\u{a}lat's Lemma observing that all the signals are bounded and square-integrable. To conclude uniform convergence (and hence, UGAS), we may appeal to  \cite[Lemma 2]{TACDELTAPE} in the following manner: firstly, we observe that all the functions are uniformly bounded (since the system is UGS and all functions). Let $c>0$ be a generic constant independent of the initial conditions and let $\rho\in \cKinfty$. Now, integrating $Y_1(t,x(t))$ we obtain that $\Norm{x_3}_2\leq c\norm{x_\circ}$. Using this and integrating $Y_2(t,x(t))$ and $Y_3(t,x(t))$ we have that $\max\{ \Norm{\phi_1(t,x(t))}_2,\, \Norm{\phi_2(t,x(t))}_2\}  \leq c\norm{x_\circ}$. In a similar way, using this we obtain by integrating $Y_4(t,x(t))$ and $Y_5(t,x(t))$ that $\Norm{x_1}_2\leq \rho(\norm{x_\circ})$ and $\Norm{x_2}_2\leq \rho(\norm{x_\circ})$. So UGAS follows observing that $\max\{\Norm{x}_\infty\,,\, \Norm{x}_2 \leq \rho(\norm{x_\circ})\}$ for all $(t,x)\in \mR\times {\cal B}(\Delta)$.
\qed

Notice that what is crucial in the argumentation above is to have that each function $Y_k(t,x)$ is nonpositive (which is less restricitve than requiring sign-definiteness) on the sets where all the functions $Y_j$ with $j<k$ vanish. Also, it is implicitly required that all the functions $Y_k(\cdot,\cdot)$ vanish only on the set $\{x = 0,\,\phi_j=0 \}$. Notice however, that $Y_k(t,x)$ are allowed to depend on time through the continuous functions $\phi_j(t,x)$ and therefore they are not required to be definitely nonzero on the sets where $Y_j$ with $j<k$ vanish. Finally, we remark that to conclude uniform convergence, all these properties are required to hold only on compact sets of the states since the system is UGS. This suggests that alternatively to the integral conditions evoked above, one may use directly the properties of the functions $V_i$ and $Y_i$. This is established in next section.

\section{Matrosov's theorem using a family of auxiliary functions}
\label{sec:matrosov}

 Now we are ready to present our main theorem. Theorem \ref{thm:main:new:mat} below  generalizes in certain directions, \cite[Proposition 2]{INTLEMMCSS} (see also \cite[Proposition 2]{INTLEMNOLCOS}) which, as clearly shown in that reference is, in its turn, an extension of the ``classical'' Matrosov theorem \cite{MAT62} which combines an auxiliary function
with a Lyapunov function that establishes UGS. See also the more recent expositions \cite[Theorem 5.5, p.58]{ROUMAW}, \cite[Theorem 2.5, p. 62]{ROUHABLAL} and \cite[Theorem 55.3]{HAH}. 
\begin{theorem}
\label{thm:main:new:mat}
Under the following assumptions the origin of the system \rref{eq:nltv} is UGAS.  \\[-8mm]
\begin{assumption}
\label{ass:mat:1}
The origin of the system \rref{eq:nltv} is UGS. \\[-8mm]
\end{assumption}
\begin{assumption}
\label{ass:mat:2}
There exist integers $j$, $m > 0$ and for each $\Delta >0$  there exist
\begin{itemize}
\item a number $\mu >0$
\item locally Lipschitz continuous functions\\ $V_i:\mR\times\mR^n\to \mR$, $i\in \{1,\ldots,j \}$ 
\item a (continuous) function $\phi:\mR\times\mR^n\to \mR^{m}$,\\ $i\in \{1,\ldots,m \}$
\item continuous functions $Y_i:\mR^n\times\mR^m\to \mR$,\\ $i\in \{1,\ldots,j \}$
\end{itemize}
such that, for almost all $(t,x) \in \mR\times \BD$, 
\begin{eqnarray}
\label{bndonViPhii}
&\dty \max\left\{ \norm{V_i(t,x)},\ \norm{\phi(t,x)} \right\}  \, \leq\, \mu, & \\
\label{bndondotVi}
&\dty \dot V_i(t,x) \,\leq\, Y_i(x,\phi(t,x)) &  \ .
\end{eqnarray}
\end{assumption}
\begin{assumption}
\label{ass:mat:3} 
For each integer $k \in \left\{1,\cdots,j\right\}$ we have that
\begin{description}
\item{(A):}\ 
$\left\{ \quad  (z,\psi) \in \BD \times \Bmu  \ , \ Y_i(z,\psi) = 0 \right. $
\ \ $ \left. \forall\, i\in \{ 1,\ldots, k-1\}  \quad \right\}$
\end{description}
implies 
\begin{description}
\item{(B):}
$\quad \{ \ Y_k(z,\psi) \leq 0 \ \}$ .
\end{description}
\end{assumption}
\begin{assumption}
\label{ass:mat:4} 
We have that
\begin{description}
\item{(A):}
$\quad \left\{ \ (z,\psi) \in \BD \times \Bmu \ , \quad
\ Y_i(z,\psi) = 0  \right. $ 
\ \ $ \left. \forall \, i\in \{ 1,\ldots, j\} \ \right\}$
\end{description}
implies
\begin{description}
\item{(B):}
$\quad \{ \ z = 0 \ \}$ .
\end{description}
\end{assumption}
\end{theorem}

In certain cases it is also useful to remark that the condition of uniformity in the initial states implicit in Assumption \ref{ass:mat:1}  can be relaxed as follows:
\begin{theorem}
\label{theorem:2}
If Assumption \ref{ass:mat:1} in Theorem \ref{thm:main:new:mat} is replaced by
\begin{quote}
\begin{enumerate}
\item
The origin is uniformly stable;
\item
for each $x_{\circ}$ there exists $M_{\circ}>0$ such that $|x(t,t_{\circ},x_{\circ})| \leq M_{\circ}$ for all $t_{\circ}$ and $t \geq t_{\circ}$;
\item
the trajectories are continuous in the initial state $x_{\circ}$ uniformly in $t_{\circ}$
\end{enumerate}
\end{quote}
and Assumptions \ref{ass:mat:2}-\ref{ass:mat:4} hold then,  the origin of \rref{eq:nltv} is UGAS.
\end{theorem}
A corollary of this result generalizes to time-varying systems
well-known results
for nonlinear cascades:
\begin{eqnarray}
\label{casc:x} \dot x &=& f(t,x,z)\\
\label{casc:z}\dot z &=& g(t,z)\,.
\end{eqnarray}
\begin{corollary}
If, for \rref{casc:x}, \rref{casc:z},
each initial condition $(x_{\circ},z_{\circ})$ produces trajectories
that are bounded uniformly in the initial time,
the functions $f$ and $g$ are locally Lipschitz uniformly in $t$, 
and the origins of \rref{casc:z} and $\dot x = f(t,x,0)$ are UGAS then, 
the origin of \rref{casc:x}, \rref{casc:z} is UGAS.
\end{corollary}

\subsection{Example: stabilization of nonholonomic systems revisited}
\label{sec:exnonhol}
With aim of building intuition for the communications channels control problem which we present in next section and of illustrating further the utility of Theorem \ref{thm:main:new:mat}, let us come back to the example of the chained-form systems \rref{eq:nonhol}. We now present and discuss  a result for the more general case of $n$ states i.e.,
\begin{subequations}\label{eq:nonhol:n}
\begin{eqnarray}
\dot x_1 &=& u(t,x)\\
\dot x_2 &=& u(t,x)x_3\\
  & \vdots &\\
\dot x_{n-1} &=& u(t,x)x_{n-1}\\
\dot x_n &=& v(t,x)
\end{eqnarray}
\end{subequations}
where $n>3$. 

We will use the following smooth control laws which are the counterparts for $n$-states, of the controller \rref{control}:
\begin{subequations}\label{control:nstates}
\begin{eqnarray}
\label{control:nstates:a} v(t,x) & =& - k_n' x_n - k_{n-1}' u(t,x)x_{n-1} - k_{n-2}'x_{n-2} -  k_{n-3}' u(t,x)x_{n-3} -\cdots \\
\label{control:nstates:b} u(t,x) &= & - k_1x_1 + h(t,z)\,, \qquad z:=\col[x_2\,,\ldots\,,x_n]\,.
\end{eqnarray}
\end{subequations}
where $k_i'>0$ for all $i\leq n$. That is, the occurrence of $u(t,x)$ in $v(t,x)$ alternates and the last term of $v(t,x)$ is $k_2' u(t,x)x_2$ if $n$ is odd or, $k_2'x_2$ if $n$ is even. Interestingly, the last $n-1$ equations of the closed loop system has the following form which recalls us of the controllability canonical form of linear systems: 
\begin{equation}\label{canonical}
\begin{bmatrix}
\dot x_2 \\ \dot x_3 \\ \vdots \\ \dot x_{n-1} \\ \dot x_n
\end{bmatrix} =
\begin{bmatrix}\,
0  &  u  & 0 & \cdots & 0 \\  0  &  0  & u & \ddots & \vdots
\\ \vdots  &  \ddots  & \ddots & \ddots &  0 \\  0  &  \cdots  & 0 & 0 & u \\   &
 \cdots  & -k_{n-2}' & -k_{n-1}'u & -k_n
\end{bmatrix}
\begin{bmatrix}
x_2 \\  x_3 \\ \vdots \\ x_{n-1} \\  x_n
\end{bmatrix}\,.
\end{equation}

The analysis problem for this case is much more complex than for systems with three states. One of the main difficulties is that the high relative degree of the system from an artificial input $w$ added to the second equation, to the output $x_n$. Notice that this relative degree in the case of 3 states is equal to 1. This is what permitted in \cite{uped}, to use certain results from linear systems theory and particularly from theory of strictly positive real systems (passive systems). Therefore, this high relative degree is a major stumbling block in the stabilization of systems with more than 3 states.

An elegant solution to this structural obstacle is the time and state-dependent coordinate transformation proposed in \cite{SAM3} which transforms the system \rref{canonical} into a system whose last $n-1$ equations have the skew-symmetric form \rref{eq:1:bis}. As remarked also in \cite{NONHOLEJC} this can also be accomplished by defining recursively, 
\begin{eqnarray}
  \label{xnu2n}
  \bar x_2 & = & \x2 \\
  \bar x_3 & = & \x3 \non \\
  v_{3} & = & -\k2u\bar x_2 - \k3\bar x_3 \non \\
  \bar x_i & = & -\jac{}{u}(v_{i-1}) + x_i, \ \ \mbox{for all } 3 < i
  \leq n  \non \\
   v_i & = & \frac{d}{dt}\left(\jac{}{u}(v_{i-1})\right) - k_{i-1}u_1\bar x_{i-1} - k_i\bar x_i, \ \ \mbox{for all } 3 < i \leq
   n   \non \\ v & = & v_n(u,\bar x)\, \label{uu2nu}
\end{eqnarray}
where the constants $\{k_2,\ldots,\,k_n\}$ are all positive and are uniquely defined, given $n$, and a  set of positive gains $\{k_2',\ldots,\,k_n'\}$. This follows from the fact that the expressions for the new variables $\bar x_i$ depend linearly in $u$ and $x$. Using the set of equations \rref{xnu2n}--\rref{uu2nu} one can write the equivalent closed loop dynamics in the ``skew symmetric'' form
\begin{subequations}\label{eq:1:bis}
\begin{eqnarray}
\label{eq:1:bis:a}
     \dot x_1 & = & -k_1 x_1 + h(t,z)=: u(t,x_1,\bar x_{2,n}) \\
\label{eq:1:bis:b}
     \dot{\bar x}_{2,n} & = & A(u(t,x_1,\bar x_{2,n})) \bar x_{2,n}
\end{eqnarray}
\end{subequations}
where  $\bar x_{2,n}:=\col[\bar x_2,\,\ldots,\,\bar x_n]$ and 
\begin{equation}\label{aofu}
   A(u) :=  \left[
            \begin{array}{ccccccc}
                0 & u & 0 & \cdots &  0 \\
               -k_2u & \ddots & \ddots & \ddots & \vdots \\
               0 & \ddots &  \ddots & \ddots & 0  \\
               \vdots & \ddots & \ddots & \ddots & u \\
                0 & \ldots &  0   & -k_{n-1}u    &  -k_n
            \end{array}
         \right] \ .
\end{equation}

To conclude UGAS for \rref{eq:1:bis} we could reason in the following manner: consider the argument of the matrix in \rref{eq:1:bis:b} along $x_1(t)$ that is, $u(t,x_1(t),\bar x_{2,n})$ then, the system \rref{eq:1:bis} may be regarded as a cascaded time-varying system. Following this line of thought it is apparent that if $h(\cdot,\cdot)$ satisfies some appropriate boundedness and regularity properties, the main question is how to guarantee that the origin of the system in \rref{eq:1:bis:b} is UGAS.

UGAS (uniformly also in $\bar x_1$) of \rref{eq:1:bis:b} may be established by a recursive coordinate transformation plus output injection arguments which transforms the system into the form 
\begin{equation}\label{lss}
\dot \xi =  \left[
            \begin{array}{ccccccc}
                -u^2 & u & 0 & \cdots  & 0 \\
               -k_2u & -u^{4} & \ddots & \ddots & \vdots \\
                0 & \ddots & \ddots & \ddots & 0 \\
                \vdots & \ddots & \ddots   & -u^{2n} & u \\
                0 & \ldots   &  0   & -k_{n-1}u    &  -k_n
            \end{array}
         \right] \xi + K(t,\xi)
\end{equation}
where $K(t,\xi)$ is uniformly square integrable along trajectories. Notice that for the particular case of\footnote{Obviously the same holds with a properly weighted Lyapunov function and $k_i\neq 1$.} $k_i=1$ we have that the derivative of $V=0.5\norm{\xi}^2$ yields 
$$
\dot V = - \sum_{i=1}^n u^{2i} \xi_i^2 + \xi^\top K(t,\xi)\,,\qquad i\leq n-1\,.
$$
Then, it can be shown that if  $u^{2i}$ is U\ped\ for any $i$, 
UGAS of \rref{lss} may be concluded owing to the fact that $K(t,\xi(t))$ is uniformly square integrable and the origin is already UGS. Finally, we see that from the properties of U\ped\ functions it follows that $u^{2i}$ is U\ped\ for any $i$ if $u$ is U\ped. See \cite{NONHOLEJC} for further details.

While the proof in \cite{NONHOLEJC} allows to conclude UGAS under very similar conditions than those imposed in \cite{SAM3}, it is rather complex since it appeals to a number of auxiliary results for time-varying systems. For instance, in the proof of the main result in \cite{NONHOLEJC} the system is regarded as a parameterized time-varying {\em linear} system $\dot \xi = A(u(t,x_1(t,\lambda),\bar x_{2,n}(t,\lambda))\xi$ where $\lambda = \col[t_\circ,x_\circ]$. The theorem that we present below can actually be proved directly using Theorem \ref{thm:main:new:mat}. 
Let us consider the system \rref{eq:nltv} with 
\begin{equation}\label{eq:41}
    x := \left[ \begin{array}{c}
                y \\
                z
               \end{array}
        \right]   \ , \qquad f(t,x) := \left[
               \begin{array}{c}
                    - y + h(t,z) \\
                    A(u(t,y,z)) z
               \end{array}
                 \right] \,, \quad z\in\mR^m\,,\ y\in\mR\,.
\end{equation}
and, imilarly to \rref{defofpsi},  define
\begin{equation}
   \label{defofpsi:nstates}
   \varphi(t,\xi) := \jac{h_\circ}{t}\,, \qquad 
h_\circ(t,\xi) := \frac{\partial h}{\partial t}\left(t,\left[
     \begin{array}{c}
           \xi \\ 0 
     \end{array}
       \right] \right)\,, \qquad 
\xi := \left[ \begin{array}{c}
                  z_{1} \\
                  \vdots \\
                  z_{m-1}
                     \end{array}
               \right]\,.
\end{equation}
Then, we have the following.
\begin{theorem}\label{thm:skew-sym}
If the function $t \mapsto \displaystyle \varphi(t,\xi)$ defined above is U\ped, the functions $z \mapsto h(t,z)$
and $\xi \mapsto \varphi(t,\xi)$ are locally Lipschitz uniformly in
$t$ and $h(t,0)=0$ then the origin of \rref{eq:nltv}, \rref{eq:41}  is UGAS. Moreover, the origin is UGAS only if the control input $u$ is U\ped\ with respect to $z$.
\end{theorem}
A direct consequence is the following.
\begin{proposition}\label{nonhol:nstates}
Consider the system \rref{eq:nonhol:n} in closed loop with \rref{control:nstates}. Let
\begin{assumption*}
the map $(t,z)\mapsto h(t,z)$ be such that $h(t,0)\equiv 0$, all its first and second partial derivatives be uniformly bounded by $\rho(\norm{z})$ where $\rho(\cdot)$ is a non decreasing function
\end{assumption*}
and let the function  $t \mapsto \displaystyle \varphi(t,\xi)$ is U\ped. Then, the closed loop system \rref{eq:1:bis} is UGAS.
\end{proposition}

The proof of Theorem \ref{thm:skew-sym} is not presented here since it follows as a particular case of the proof of Theorem \ref{thm:channels} presented in next section along the lines of the proof for the case of three states. However, for completeness and further reference we provide the guidelines for sufficiency. This follows by applying directly Theorem \ref{thm:main:new:mat} with the following functions.
\begin{equation}
  V_{1}(t,x) := z^\top P z\,
\end{equation}
where $P:=\diag\{p_1,\,\cdots,\,p_m\}$ with $p_i:=\dty\frac{p_{i-1}}{k_{i}}$ for all $ i \in [2,\,\ldots,\,m-1]$, $p_m = p_1 := 1$. In particular, this function allows to show UGS. The rest of the auxiliary functions are:
\begin{eqnarray}\label{vs:nonhol}
   V_{i}(t,x) & := & z_{m-i+2} \cdot u(t,y,z)^{2i-3} \cdot z_{m-i+1}, \qquad i=2,\cdots,m\\
   V_{m+1}(t,x) & := & \zeta(t,x)^{2}
\end{eqnarray}
with 
\begin{eqnarray*}
  \omega(t,\xi) &:=&  \int_{-\infty}^{t} e^{-(t-\tau)} \varphi(\tau,\xi) d\tau\\
   \zeta(t,x) &:=& y - \widetilde{h}(t,\xi)   + \omega(t,\xi)\\
  \widetilde{h}(t,\xi) &:=&
h\left(t,\left[ \begin{array}{c}
          \xi \\
           0
           \end{array}
            \right] \right)\,.
\end{eqnarray*}
Defining for $i=2 \cdots m$, $\widetilde{\phi}_{i}(t,\xi) := \omega(t,\xi)^{i-1} z_{m-i+1}$,
\begin{eqnarray}\label{eco}
  V_{m+i}(t,\xi) &:=&
   - \int_{t}^{\infty} e^{t-\tau} |\widetilde{\phi}_{i}(\tau,\xi)|^{2} d\tau
\end{eqnarray}
and finally, 
\begin{equation}
  V_{2m+1}(t,x) := y^{2}\,.
\end{equation}
Note the similarity with the functions for the case of three states. In particular, note that here we need $n-2$ functions as defined in \rref{vs:nonhol} and $n-2$ functions as in \rref{eco} (see also \rref{v2} and \rref{v4}\,) where $n-2$ corresponds to the relative degree of \rref{eq:1:bis:b}.

\section{Control of communication channels}
\label{sec:channels}
\def\figpath{.}

We will address the problem of stabilizing by smooth feedback, a series of port-interconnected driftless systems as illustrated in  Figure \ref{channels}. The setting of this control problem covers that of time-varying smooth feedback stabilization of chained-form nonholonomic systems (cf. \cite{KOLMCC,MURSAS}) and we will solve it by appealing to Theorem \ref{thm:main:new:mat}.

\subsection{Problem setting and its solution}

 \putfigtex{channels}{channels}{Communication channels.}

In Figure \ref{channels}, each block contains a bank of integrators with nonlinearities at the input and output of the integrator.  The nonlinearities are supposed to be such that each block is passive in a specific sense.  In particular, we assume that each block can be modeled as
\begin{equation}
   \Sigma_{i} : \qquad  \left\{ 
  \renewcommand{\arraystretch}{1.6}
  \begin{array}{rcl}
   \dot{x}_{i} & = &  B_{i}(x_{i}) u_{i} \\
       y_{i} & = & h_{i}(x_{i}) =  B_{i}(x_{i})^{T} \nabla W_{i}(x_{i}) \\
  \end{array}
   \right.
\label{eq:block}
\end{equation}
where $i\in \{1,\ldots, n\}$ $y_i$, $u_i\in\mR^{p_i}$, $x_i\in\mR^{n_i}$, $B_i(\cdot)$ and  $\nabla W_{i}(\cdot)$ are locally Lipschitz and the following bounds hold:
\begin{eqnarray}
\label{nablahW}
   \nabla h(x_{i}) B_{i}(x_{i}) & \geq & c_{i} I\,,   \qquad \qquad\ c_{i} > 0 \\
\label{Bbnd}
   \norm{B_{i}(x_{i})} & \leq  & \rho_{B_i}(\norm{x_i})\,,\qquad   \rho_{B_i} \in {\cal N}  \\
\label{Walpha}
  W_{i}(x_{i}) & \geq & \alpha_{i}(|x_{i}|)\,, \qquad \ \alpha_{i} \in \cKinfty  \\
\label{rhoW}
  \nabla W_{i}(x_{i}) & \leq &  \rho_{W_i}(|x_{i}|)\,, \qquad  \rho_{W_i} \in {\mathcal K}_{\infty}  \\
\label{ykappa}
\norm{h_i(x_i)} & \geq & \kappa_{i}(|x_{i}|)\,, \qquad\ \kappa_{i} \in {\mathcal P}{\mathcal D}\,.
\end{eqnarray}
The nonlinear integrator blocks are interconnected via static, nonlinear, time-varying ``communication channels''.  In particular, the connection between  blocks $i$ and $i+1$ is modeled by the nonlinear gain function $g_{i}$ which takes values in $\mR$ and may depend, in general, on any of the states, time, and perhaps some additional states from outside  of the network.  This means that the input from the left to the $i$th block, denoted $u_{i,\ell}$, and the input from the right to the $i$th block, denoted $u_{i,r}$,
are 
\begin{eqnarray*}
   u_{i,\ell} & = & g_{i-1} y_{i-1} \\
   u_{i,r} & = & g_{i}  y_{i+1}\,.
\end{eqnarray*}
The blocks are such that the input to the nonlinear integrator is given by
\[
  u_{i} = u_{i,r} - u_{i,\ell} \ .
\]
We assume that the communication channel gains have the functional form
\begin{subequations}  \label{eq:fform}
\begin{eqnarray}
  \label{eq:fform:a}
   \dot{z}_{i} & = & -z_{i} + \widetilde{g}_{i,a}(t,x)   \\
  \label{eq:fform:b}
  \widetilde{g}_{i}(t,x,z) & = & - z_{i} + \widetilde{g}_{i,a}(t,x) + \widetilde{g}_{i,b}(t,x) \\
  \label{eq:fform:c}
    g_{i}(t,x,z) & = & \displaystyle \prod_{j=i}^{n-1} \widetilde{g}_{j}(t,x,z)
\end{eqnarray}
\end{subequations}
for $i \leq n-1$ and where the functions $\widetilde{g}_{i,a}$ 
 and  $\widetilde{g}_{i,b}$ 
are continuous and Lipschitz in $x$ uniformly in $t$.
 One interesting situation arises when all of the communication
 channel gains are the same, which is the case when
 $\widetilde{g}_{i,a}\equiv 0$, $z_{i}(0)=0$, $\widetilde{g}_{i,b}=1$
 for $i=1,\cdots,n-2$.  In this case, every gain is given
 by $\widetilde{g}_{n-1}$.

The control problem is to attach a system $\Sigma_{n+1}$ to the right of $\Sigma_{n}$, and give necessary and sufficient conditions on the communication channel gains to guarantee that the origin for the closed-loop system is uniformly globally asymptotically stable.  
For the controller, we will use any static strict ``first and third sector'' nonlinearity $\sigma(\cdot)$, and the connection to $\Sigma_{n}$ will be made with a reliable communication channel, e.g., $g_{n} \equiv 1$. In particular, we have
\begin{equation}
  \Sigma_{n+1} : \quad  \left\{ 
  \renewcommand{\arraystretch}{1.6}
  \begin{array}{ccl}
          y_{n+1} & = & \sigma(u_{n+1}),\qquad\ \, s>0 \Rightarrow \sigma(s)s > 0,\, \sigma(0) = 0 \\
         u_{n+1}^{T} \sigma(u_{n+1}) & \geq & \rho(|u_{n+1}|), \qquad
                            \rho \in {\mathcal P}{\mathcal D} \\
          u_{n+1} & = & y_{n} \\
          y_{n+1} & = & u_{n} \ .
  \end{array}
   \right.
\end{equation}

With this controller architecture and functional form for
the communication channel gains indicated in Figure \ref{channels},
 we ask the question:
\begin{quote}
{\em What are necessary and
sufficient conditions on the communication channel gains to guarantee
uniform asymptotic stability of the origin for the system
(\ref{eq:block})-(\ref{eq:fform})?}
\end{quote}

The answer will be expressed in terms of the notion of uniform $\delta$-persistency of excitation (cf. Definition \ref{def:uped}).  
Note that in  Definition \ref{def:uped}, $x$ is a constant parameter hence, the necessary and sufficient conditions for stability will be expressed in terms of the state $x$, being constant. To that end note that when $x$ is constant the $z_{i}$ subsystems in (\ref{eq:fform}) are time-invariant linear systems with time-varying inputs. To better see this, let $x= \bar x$ with $\bar x$ a constant vector and call $\bar z$ the new state of the {\em linear} system (\ref{eq:fform}) in this setting. Then, it is direct to see that the $i$th communication channel gain in (\ref{eq:fform:b}) satisfies (when $\widetilde g_{i,b}(t,x) \equiv 0$\, ) 
\begin{eqnarray}
\dot{\tilde g}_i(t,\bar x,\bar z) & = &  -\bar z_i + \tilde g_{i,a}(t,\bar x) + \dot{\wt g}_{i,a}\\
\label{filteredg}
 & = &  -{\tilde g}_i(t,\bar x,\bar z) +  \jac{\tilde g_{i,a}}{t}
\end{eqnarray}
and we note that the steady-state solution of \rref{filteredg} is given by 
\begin{equation}\label{ch:omegadef}
\omega_i(t,\bar x):= \int_{-\infty}^t \mbox{e}^{-(t-\tau)} \psi_i(\tau,\bar x) d\tau\qquad i \in[1,\ldots,n-1]
\end{equation}
where 
\begin{equation}\label{ch:psiidef}
\psi_i(t,\bar x):=  \frac{\partial \widetilde{g}_{i,a}(t,\bar x)}{\partial t}\,.
\end{equation}
That is, $\omega_i(t,\bar x)$ is the steady-state value of the $i$th communication channel gain and correspondingly, the steady-state value of the gain of the first communication channel can be computed to be the product of all the $\omega_i(t,\bar x)$'s for all $i\leq n-1$. Based on these observations and the extended Matrosov's Theorem \ref{thm:main:new:mat}\ we may establish the following.
\begin{theorem}\label{thm:channels}
Suppose that 
\begin{assumption*}\label{assongiab}
the   function $\wt g_{i,b}(\cdot,\cdot)$ is continuous, locally Lipschitz in $x$ uniformly in $t$, and bounded uniformly in $t$. The function $\wt g_{i,a}(\cdot,\cdot)$ is locally Lipschitz in $x$ uniformly in $t$, twice continuously differentiable with first and second partial derivatives locally Lipschitz in $x$ uniformly in $t$ and $\tilde g_{i,a}(t,0) \equiv 0$.
\end{assumption*}
Then, the origin of the system (\ref{eq:block})-(\ref{eq:fform}) is uniformly globally asymptotically stable if and only if the functions
\begin{equation}\label{eq:thm:main}
(t,x_{1},\ldots,x_i) \mapsto   
\left. \prod_{j=i}^{n-1} 
 \left(\omega_j(t,x)  + \widetilde{g}_{j,b}(t,x)  \right)
 \right|_{x_{i+1}=\ldots=x_{n}=0}  \qquad \forall i \in
 \left\{1, \ldots, n -1 \right\}
\end{equation}
is $U\delta$-PE with respect to $x_i$. 
\end{theorem}
Roughly speaking, the sufficient condition for UGAS is that the first  communication channel gain be U$\delta$-PE, at least when the last ${n-1}$ states are zero.  From the structure of the gains in (\ref{eq:fform}) and Fact \ref{fact:prodphis}, this implies that {\em each} of the communication channel gains is U$\delta$-PE when the last ${n-1}$ states are zero. This does not mean that it is required that each channel is always functioning but rather that each communication channel is functioning on average and in a synchronized way.  Moreover, this average should be uniform in time.  However, the average does not need to be uniform in the state.  For example, the quality of the communication channel could possibly degrade as the ``power'' of the transmission signals, perhaps encoded by the size of the states $x_{i}$, decreases to zero. This idea is captured by the notion of U\ped, in that this property roughly means for a function that it is PE in the usual sense (the channel work in average) but with a ``degree of excitation'' which depends on the size of the state (cf. \cite[Lemma 3]{DPECDC02}.

\begin{remark}
It is worth mentioning  that the system architecture above covers the so-called ``skew-symmetric'' systems considered in \cite{SAM3,NONHOLEJC} and in Section \ref{sec:exnonhol}. We may see this if we let $\tilde g_i = u_1$ for all $i$ where $u_1$ is one of the two control inputs in those references (in particular, the controller which is required to be U\ped\ in \cite{NONHOLEJC}), $y_i=x_i$ for each $2\leq i \leq n$ we replace $z_i$ by\footnote{ Here, the index $i$ does not make sense because we have only one $z$-state, this is because all the gains $\tilde g_i = u_1$.} $ x_1$ and finally, we relate the function $\tilde g_{i,a}$  to the function whose second derivative in \cite{NONHOLEJC} is required to be U\ped\ or, to the ``heat function'' in \cite{SAM3}.
\end{remark}

\subsection{Proof of Theorem \ref{thm:channels}}

We analyze the stability of the origin of \rref{eq:block}, \rref{eq:fform} with state $(x,z)$ by analyzing that of the following equivalent system. We consider the dynamics of \rref{eq:block} and, instead of the channel gain dynamics, \rref{eq:fform}, we consider that of the difference between the gains and their steady state solution (without the additional $g_{i,b}$ term). More precisely, let 
\begin{equation}
\label{ch:zeta}
\zeta_i(t,x,z):= z_i - \widetilde g_{i,a}(t,x) + \omega_i(t,x)\,, \qquad i \in[1,\ldots,n-1]
\end{equation}
then (droping the arguments for simplicity in the notation),
\begin{equation}\label{ch:dotzeta}
\dot \zeta_i = - z_i + \widetilde g_{i,a} - \jac{ \widetilde g_{i,a}}{t} - \jac{ \widetilde g_{i,a}}{x}\dot x + \jac{\omega_i}{t} +  \jac{\omega_i}{x}\dot x
\end{equation}
where  $\dot x = \col[\,B_i(g_iy_{i+1}-g_{i-1}y_{i-1})\,] $ and 
\begin{eqnarray}\label{qqq1}
- \jac{ \widetilde g_{i,a}(t,x)}{t} & = & -\psi_i(t,x)\\
\jac{\omega_i(t,x)}{x} & = & \int_{-\infty}^t \mbox{e}^{-(t-\tau)} \jac{\psi_i}{x}(\tau,x) d\tau \\
\label{qqq2}
\jac{\omega_i(t,x)}{t} & = & \psi_i(t,x) - \omega_i(t,x)\,.
\end{eqnarray}
Hence, using \rref{ch:zeta}, \rref{qqq1}-\rref{qqq2}  in \rref{ch:dotzeta} we obtain that
\begin{equation}
\label{ch:dotzeta2}
\dot \zeta_i = -\zeta_i  + \left[\int_{-\infty}^t \mbox{e}^{-(t-\tau)} \jac{\psi_i}{x}(\tau,x) d\tau - \jac{ \widetilde g_{i,a}}{x} \right] \dot x\,, \qquad i \in \{1,\ldots, n-1\}\,.
\end{equation}
Under Assumption \ref{assongiab} we have that the right hand side of \rref{ch:dotzeta2} is  continuous in $t$ and locally Lipschitz in $(\zeta,x)$ uniformly in $t$. From the fact that $\wt g_{i,a}(t,0)\equiv 0$ and $\wt g_{i,a}(\cdot,x)$ is continuously differentiable we also have that $\psi_i(t,0)\equiv 0$ and in view of \rref{ch:omegadef}, so does $\omega_i(t,0)\equiv 0$. Hence, it follows from \rref{ch:zeta} that  $z_i=0$ if $\zeta_i=0$ and $x_i=0$. 

It is left to show that the necessary and sufficient condition for the origin $(\zeta,x)=(0,0)$ of \rref{eq:block}, \rref{ch:dotzeta2} to be UGAS is that the functions in \rref{eq:thm:main} be U\ped\ with respect to $x_i$ for all $i\in \{1,\ldots, n-1\}$.

\subsubsection{Necessity}

We invoke the condition that for the origin of a system $\dot{\xi}=F(t,\xi)$ with $F(\cdot,\cdot)$ locally Lipschitz in $x$ uniformly in $t$ and continuous in $t$ to be UGAS, necessarily $F(\cdot,\xi)$ must be $U\delta$-PE with respect to $\xi$ (cf. \cite[Theorem 2]{DPEALAMOS}).  Let $F(t,\xi)$ correspond to the right hand side of \rref{eq:block}, \rref{ch:dotzeta2} with 
\begin{equation}
   \xi:= \left[ \begin{array}{c}
               \zeta \\
                  x
                \end{array}
         \right] \ 
\end{equation}
and consider the partition $F(t,x):= \col[ F_\zeta(t,\xi);\,F_x(t,\xi)]$ where $F_\zeta(t,\xi)$ corresponds to the vector whose elements correspond to the right hand side of \rref{ch:dotzeta2} with $i \in \left\{1, \ldots, n-1\right\}$. Let $i$ be arbitrary in such interval, and let us consider nonzero points $\xi$ such that $\zeta=0$ and $x_{i+1}=\ldots=x_{n}=0$. For ease of reference let us define these points as $\xi=s$ and define the following functions
\begin{eqnarray}
\beta_i(t,x) & :=  & \left[\int_{-\infty}^t \mbox{e}^{-(t-\tau)} \jac{\psi_i}{x}(\tau,x) d\tau - \jac{ \widetilde g_{i,a}}{x} \right]\,. 
\end{eqnarray}
Hence, letting $(\cdot)^\circ$ denote $(\cdot)$ evaluated at $s$, we have that 
\begin{equation}
\label{Fofts}
  \begin{bmatrix}
    \dot \zeta_1^\circ \\[1mm] \vdots\\[1mm] \dot \zeta_n^\circ  \\[1mm]
    \dot x_{1}^\circ  \\[1mm] 
    \dot x_{2}^\circ  \\[1mm] \vdots\\[1mm]
    \dot x_{i-1}^\circ  \\[1mm]
\dot x_{i}^\circ  \\[1mm]
\dot x_{i+1}^\circ  \\[1mm]
\dot x_{i+2}^\circ  \\[1mm]
\vdots\\
\dot x_n^\circ 
  \end{bmatrix}
= \begin{bmatrix}
    \beta_1(t,s)  F_x(t,s)\\
    \vdots\\
    \beta_n(t,s)  F_x(t,s)\\[1mm]
    B_1(x_1)[\,\widetilde g_{1}^\circ \cdots \widetilde g_{n-1}^\circ \,]h_2(x_2)\\[1mm]
    B_2(x_2)[\,-\widetilde g_{1}^\circ \cdots \widetilde g_{n-1}^\circ h_1(x_1) + \widetilde g_{2}^\circ \cdots \widetilde g_{n-1}^\circ h_3(x_3)\, ] \\[1mm]
    \vdots\\[1mm]
    B_{i-1}(x_{i-1})[\,-\wt g_{i-2}^\circ\cdots\wt g_{n-1}^\circ h_{i-2}(x_{i-2})  + \wt g_{i-1}^\circ\cdots\wt g_{n-1}^\circ h_{i}(x_{i}) \,] \\[1mm]
    B_{i}(x_{i})[\,-\wt g_{i-1}^\circ\cdots\wt g_{n-1}^\circ  \,] h_{i-1}(x_{i-1})\\[1mm]
    B_{i+1}(0)[\,-\wt g_{i}^\circ \cdots\wt g_{n-1}^\circ \,] h_i(x_{i}) \\[1mm]
  0\\ \vdots \\ 0\,
  \end{bmatrix}  
\begin{array}{c}
\ \\[1.4cm]
  \shiftleft{9mm}
\left.
\begin{array}{c}
  \ \\[5cm]
\end{array}
\right\}= F_x(t,s)\,.
\end{array}
\end{equation}
One can see from \rref{Fofts} that the function
\begin{equation}
      \left.      
\left( \prod_{j=i}^{n-1} \wt{g}_{j} \right) 
\right|_{x_{i+1}=\ldots=x_{n}=0} \quad \mbox{ for any } i\in\{1,\ldots,n-1\}
\label{eq:5}
\end{equation}
can be factored out of $F(t,s)$. Then, observing that for any real number  $\alpha$ and any real vector $v$ we have that $\norm{\alpha v}\leq\norm{\alpha}\norm{v}$ and that $F(t,s)$ is U\ped\ if and only if so is $\norm{F(t,s)}$, we invoke Fact \ref{fact:prodphis} to obtain that necessarily any of the functions in \rref{eq:5} is U\ped\ with respect to $x_i$. The result follows using \rref{ch:zeta} and \rref{eq:fform:b} to see that $\zeta_j = 0 $ implies that $\wt g_j = \omega_j + \wt g_{j,b}$.

\subsubsection{Sufficiency}

\noindent \underline{Proof of UGS:}
  The Lyapunov function
\begin{equation}
   V_{1}(x) := \sum_{i=1}^{n} W_{i}(x_{i}) \ ,
\end{equation}
which is positive definite and proper with respect to $x$,
has the property that
\begin{equation}
   \dot{V}_{1}(x) \leq - \rho(|y_{n}|) \leq 0
\label{eq:Lyap1}
\end{equation}
regardless of the properties of the communication
channels.  So, with the guarantee of local existence of solutions, 
we have that
\begin{equation}
   |x(t)| \leq  \gamma(|x_{\circ}|)  \qquad \forall t \geq t_{\circ} \ , \
               x_{\circ} \in \mR^{n} \ .
\end{equation}
Technically, we only have this bound on the maximal interval of
definition.  But with $x$ bounded on the maximal interval of definition
and Assumption \ref{assongiab}, it follows that
\begin{equation}
  |\zeta_{i}(t)| \leq |\zeta_{i}(t_\circ)| + \gamma(||x||_{\infty})
\end{equation}
on the maximal interval of definition.  Thus, solutions are
defined for all time and, in fact, the origin is UGS. 

\begin{remark}
  Notice that the same conclusion follows for the origin $(z,x)=(0,0)$ following similar arguments and using \rref{eq:fform:a}.
\end{remark}

\noindent \underline{Proof of UGA} (with a family of sign-indefinite Lyapunov functions): we will define $2n-2$ functions which may be classified in 3 groups of $n-1$ functions. The first group of functions is defined as follows\footnote{Strictly speaking we should write $V_2(t,x,\zeta) := h_n(x_n) \cdot g_{n-1} (\,t,x,\zeta - \omega(t,x) + \wt g_{a}(t,x)\,) \cdot h_{n-1}(x_{n-1})$ and similarly for the rest of the functions but we prefer to drop the arguments and make the obvious abreviations to avoid a cumbersome notation.}
\begin{eqnarray}
 V_{2} &:=&  y_{n} \cdot  g_{n-1} \cdot  y_{n-1} \\
 V_{3} &:=&  y_{n-1} \cdot  g_{n-2}\,g_{n-1}^2 \cdot  y_{n-2} \
\end{eqnarray}
where `$ \cdot $' denotes the scalar product and for $i\in \{4,\cdots,n\}$\,,
\begin{equation}
 V_{i} :=  y_{n-i+2}\cdot  g_{n-i+1} \,g_{n-i+2}^{2}
   \cdots g_{n-1}^{2}\cdot y_{n-i+1} \,. 
\end{equation}
For clarity we recall that $y_{n-i+2}\in\mR^{p_i}$ and the gains $g_{k}$ are scalar functions. We proceed to compute some bounds for the time derivative of $V_i$. To that end we first notice that 
\begin{eqnarray}
  \label{doty}
  \dot{y}_{j} & = & \nabla h_{j}(x_{j}) B_j(x_{j})
 \left[g_{j} y_{j+1}- g_{j-1} y_{j-1} \right], \qquad j\in [2,\ldots, n]\,.
\end{eqnarray}
Long but straightforward calculations, which involve the use of the bounds \rref{nablahW}--\rref{ykappa}, show that for each $\Delta>0$ there exists  $\nu>0$ such that for all $(x,\zeta)\in {\cal B}(\Delta)^2$, we have
\begin{subequations}\label{dotvis}
\begin{eqnarray}
\label{dotvis:a} \dot{V}_{2} & \leq & - c_{n} |g_{n-1} y_{n-1}|^{2} + \nu (|y_{n}| + \norm{\sigma(y_n)} )\\
\label{dotvis:b}
 \dot{V}_{3} & \leq & - c_{n-1} |g_{n-2} g_{n-1} y_{n-2}|^{2}
                         + \nu \left( |y_{n}| + |g_{n-1} y_{n-1}| \right) \,\\
  \dot{V}_{i} & \leq & - c_{n-i+2} |g_{n-i+1} \cdots g_{n-1} y_{n-i+1}|^{2} \nonumber \\ &&  
\label{dotvis:c}
          + \nu \left( |g_{n-i+2} \cdots g_{n-1} y_{n-i+2} | + | 
       g_{n-i+3} \cdots g_{n-1} y_{n-i+3} |  \right)\,
\end{eqnarray}\end{subequations}
for all $i\in [4,\ldots,n]$ and where the coefficients $c_i$ come from \rref{nablahW}. Then, using \rref{2001} and \rref{2002} successively we obtain that for all $i \in [2\ldots n-2]$,
\begin{equation}
\norm{g_i y_{i+1}} \leq \norm{[g_{n-1}\cdots g_{i+1}]}^{\ \! 1/n-i-1}\, \norm{y_{i+1}}^{1/n-i-1} \ \norm{y_{i+1}}^{n-i-2/n-i-1} \ \norm{\widetilde g_{n-2}^{\ \! 1/n-i-1} \,\cdots\, \widetilde g_{i+1}^{\ \! n-i-2/n-i-1}\widetilde g_i }
\end{equation}
so  defining the following $n-1$ functions
\begin{subequations}\label{ch:phis}
\begin{eqnarray}
 \phi_{n-1} & := & \norm{g_{n-1}}\norm{y_{n-1}}\\
  \phi_{n-i+1} & := & \norm{g_{n-i+1} \cdots g_{n-1}}\norm{ y_{n-i+1} }\,,\qquad i\in [3,\ldots,n]
\end{eqnarray}\end{subequations}
we can rewrite the inequalities in \rref{dotvis} as
\begin{subequations}\label{dotvis2}
\begin{eqnarray}
\label{dotvis2:a} \dot{V}_{2} & \leq & - c_{n} \phi_{n-1}^{2} + \nu (|y_{n}| + \norm{\sigma(y_n)}\,) \\
\label{dotvis2:b}
 \dot{V}_{3} & \leq & - c_{n-1} \phi_{n-2}^{2}
                         + \nu \left( |y_{n}| + \phi_{n-1}\, \right) \,\\
  \dot{V}_{i} & \leq & - c_{n-i+2} \phi_{n-i+1}^{2}
          + \nu \left( \phi_{n-i+2}  +  \phi_{n-i+3}  \, \right)\,.
\end{eqnarray}\end{subequations}

\begin{remark}
We wish to emphasize  the way the functions $V_1$ to $V_{n}$ defined so far, are ordered.  For the use of Theorem \ref{thm:main:new:mat} what is important to observe is that $\dot V_2 \leq 0$ on the set where the bound on $\dot V_1$ is zero, that is when $y_n \equiv 0$. Accordingly, each of the bounds on the succeeding $\dot V_i$'s contain three essential terms: the first is a negative term of $\phi_{n-i+1}$, the next two correspond to a number $\nu$ times  $\phi_{n-i+2}$ and $\phi_{n-i+3}$ which appear squared and with sign `-'  in the previous two derivatives $\dot V_{i-1}$ and $\dot V_{i-2}$ respectively. This hints at the idea that one should be able to recursively show that if $y_n \to 0$ then, so does $\phi_{n-1}$ hence also $\phi_{n-2}$, etc.\footnote{This reasoning is similar to the arguments employed in \cite{SAM3} to prove (non uniform) convergence for skew-symmetric systems.} Later, the U\ped\ assumption on the gains will be used essentially to imply that if $\phi_i\to 0$, necessarily $x_i\to 0$.
\end{remark}

The next  group of $n-1$ functions that we introduce helps to conclude on the behavior of the gains; more precisely to show that the communication channel gains converge to their steady state solution which was briefly discussed above. To that end, let
\begin{equation}\label{ch:vzeta}
V_{n+i} := \zeta_i^2\,,\qquad i\in \{1,\ldots,n-1\}\,.
\end{equation}
The aim is to bound (on compact sets of the states $x$ and $z$) the total derivatives of these functions (i.e., along the trajectories of \rref{ch:dotzeta2}) with terms of the type $-\zeta_i^2$ plus terms involving the $\norm{\phi_i}$'s defined in the previous group of functions and which appear with sign `$-$' in the previous bounds on $\dot V_i$. 

To that end, we first observe that the terms in brackets multiplying $\dot x$ in \rref{ch:dotzeta2} are uniformly bounded in $t$ by a continuous function of the norm of the state and hence, it is bounded by a number $\nu >0$ for all $x \in {\cal B}( \Delta)$. We also find it convenient to recall that
\begin{eqnarray*}
\dot x_1  & = & B_1g_1y_2\\
\dot x_j & = & B_j \left[g_{j} y_{j+1}- g_{j-1} y_{j-1} \right], \qquad j\in [2,\ldots, n-2]\\
\dot x_{n-1} & = & B_{n-1}[g_{n-1} y_n - g_{n-2}y_{n-2}]\\
\dot x_n & = & B_n[g_n \sigma(y_n) - g_{n-1}y_{n-1}]\,.
\end{eqnarray*}
so we  use \rref{2001}, \rref{2000} and \rref{eq:block} and the fact that all the gains $g_i$ and $\wtg{}{i}$ and the functions $B_i$ are uniformly bounded in $t$ on compact sets ${\cal B}(\Delta)$ and proceed as we did for \rref{doty} to find that for all $(x,\zeta)\in {\cal B}(\Delta)^2$ and all $i\in [2,\cdots,n-2]$,
\begin{equation}
\label{bndondotx}
\norm{\dot x} \leq \nu\left[\,\sum_{i=2}^{n-2}\left(\phi_{i+1}^{1/n-i-1} + \phi_{i-1}^{1/n-i+1}\right)   + \phi_2^{1/n-2}   + \norm{y_n} +\norm{\sigma(y_n)}  + \phi_{n-2}^{1/2} + \phi_{n-1} \, \right]\,.
\end{equation}
So we can use now this inequality together with Assumption \ref{assongiab} to find a bound on the total time derivative  of the function defined in \rref{ch:vzeta} along the trajectories of \rref{ch:dotzeta2}. We obtain that for all $(t,x,\zeta)\in \mR\times {\cal B}(\Delta)^2$,
\begin{equation}
\label{ch:dotvzeta1}
\dot V_{n+i} \leq - \zeta_i^2 + \nu\left( \norm{\sigma(y_n)} + \norm{y_n}  + \sum_{j=1}^{n-1} \phi_{j}^{1/n-j} \right) 
\end{equation}
for all $i\in [1,\ldots, n-1]$ which is what we were seeking for. Notice that the positive terms above appear with negative signs in the bounds on the derivatives of the previous group of auxiliary functions.

Roughly speaking, from the following group of functions we will be able to get terms that allow us to conlude on the convergence of the states to zero provided that the $\phi_i$'s converge to zero and that the gains are persistently exciting in the specific way we imposed. For the purposes of using Theorem \ref{thm:main:new:mat}, we look for a group of functions with derivatives bounded by negative terms of $x_i$ and possibly, positive terms of $\norm{\phi_i}$ and $\norm{\zeta_i}$ which appear with negative sign in teh previous bounds. 

With this in mind, we define similarly to the function $V_4$ in \rref{31},
\begin{eqnarray}
\label{vks}
V_{2n-1+i}(t,x) & := & -\int_t^\infty \mbox{e}^{t-\tau}\norm{\Omega_{n-i}(\tau,x)h_{n-i}(x_{n-i})}^2 d\tau \\
\label{tphi}
\Omega_{i}(t,x) & := & \dty  \left. \left( \prod_{j=i}^{n-1} \omega_j(t,x) + \wt g_{j,b}(t,x)\right) \right|_{x_{i+1}=\ldots=x_{n}=0}    \,, \quad i\in [1,\ldots, n-1]\ .
\end{eqnarray}
Under Assumption \ref{assongiab} these functions are locally Lipschitz in $x$ uniformly in $t$. Furthermore, using Fact \ref{fact:onpeofpsi} and the inequality $\int_t^{t+T}\norm{f(\tau)}d\tau \leq (T\int_t^{t+T}\norm{f(\tau)}^2d\tau )^{1/2}$ we obtain that for each $\Omega_i(t,x)$  which is U$\delta$-PE with respect to $x_i$, there exists $\gamma_{i}\in {\cal K}$ such that
\begin{equation}
\label{vksbnd}
 V_{2n-1+i}(t,x) \leq -\gamma_{n-i}(\norm{x_{n-i}})\norm{h_{n-i}(x_{n-i})}^2\,. 
\end{equation}
As a matter of fact, $\gamma_i(s):= \min\left\{ s,\, \frac{e^{- \theta_{i\Delta}(s)}}{\theta_{i\Delta}(s)}\gamma_{i\Delta}(s)^2\right\}$.

On the other hand, the total time derivative of $V_{2n-1+i}(t,x)$ yields at the points of existence, 
\begin{equation}
\label{dotvks1}
\dot V_{2n-1+i}(t,x) = V_{2n-1+i}(t,x) + \norm{\Omega_{n-i}(\tau,x)h_{n-i}(x_{n-i})}^2 + \jac{V_{2n-1+i}}{x} \dot x \, \quad \mbox{a.e.}
\end{equation}
To upperbound the last term on the right hand side of this expression on compact sets of the states, we use \rref{bndondotx} and the local Lipschitz property of $V_{2n-1+i}(t,x)$ and Assumption \ref{assongiab} to bound the gradient (on its points of definition).  

For the term $\norm{\Omega_{n-i}(\tau,x)h_{n-i}(x_{n-i})}^2$ in \rref{dotvks1} we derive a bound (as before on compact sets of the states) involving $\norm{\phi_j}$, $\norm{\zeta_j}$ with $j\leq n-1$  and $\norm{y_n}$ which appear with negative sign in the bounds on the previous functions' derivatives. To that end, we use \rref{ch:zeta} and \rref{eq:fform:b} to obtain a more convenient expression for $\omega_i$ and substitute it in \rref{tphi} to see that
\[
\norm{\Omega_i }^2 = \prod_{j=i}^{n-1}\ [\, \wtg{}{j} + \zeta_j   \,]^2\,
\]
so that for all $(t,x,\zeta)\in \mR\times {\cal B}(\Delta)^2$,
$
[\, \wtg{}{j} + \zeta_j \,]^2 \leq \wtg{2}{j} + \nu \norm{\zeta_j} \,.
$
 Furthermore, using once more the uniform boundedness in $t$ of $\wtg{}{j}(t,\cdot,\cdot)$ and arguing as above, we have that for any integer $j\in [i,\ldots, n-1]$, 
\[
[\, \wtg{}{j} + \zeta_j  \,]^2[\, \wtg{}{j+1} + \zeta_{j+1}   \,]^2 \leq \wtg{2}{j}\wtg{2}{j+1} + \nu(\norm{\zeta_j} + \norm{\zeta_{j+1}} )\,.
\]
Using the fact that  $\norm{y_i} \leq \rho_{B_i}(\norm{x_i})\rho_{W_i}(\norm{x_i})$ (see \rref{Bbnd} and \rref{rhoW}), it follows that for all $(t,x,\zeta)\in \mR\times {\cal B}(\Delta)^2$, 
\[
\norm{\Omega_i y_i}^2 \leq \prod_{j=i}^{n-1}\, \wtg{2}{j}y_i^2 + \nu\sum_{j=i}^{n-1} \norm{\zeta_j} \,.
\]
Notice that in view of \rref{eq:fform:c} the bound above is exactly the same as 
\[
\norm{\Omega_i y_i}^2 \leq  g^{2}_{i}y_i^2 + \nu\sum_{j=i}^{n-1} \norm{\zeta_j} 
\]
so using the identities \rref{2002} involving $g_{i}\,y_i$ we finally obtain that  for all $(t,x,\zeta)\in \mR\times {\cal B}(\Delta)^2$, and all $i\in [1,\ldots, n-1]$\,
\begin{equation}
\label{bndontphi2}
\norm{\Omega_i y_i}^2 \leq  \nu \left(\, \phi_i^{2/n-i} + \sum_{j=i}^{n-1} \norm{\zeta_j}  \,\right) 
\end{equation}
Summarizing,  we have from \rref{dotvks1}, \rref{vksbnd} and \rref{bndontphi2} that for almost all  $(t,x,\zeta)\in \mR\times {\cal B}(\Delta)^2$ and all $i\in [1,\ldots, n-1]$,
\begin{equation}
\label{dotvksfinal}
\dot V_{2n+i-1} \leq -\gamma_{n-i}(\norm{x_{n-i}})\norm{h_{n-i}(x_{n-i})}^2 + \nu\left(\,  \phi_{n-i}^{2/i} + \norm{y_n} + \norm{\sigma(y_n)} + \sum_{\ell=1}^{n-1} \phi_\ell^{1/n-\ell} + \norm{\zeta_\ell}\, \right)
\end{equation}
where the second, third and fourth term in brackets come from bounding the last  term of \rref{dotvks1} like we did to obtain \rref{ch:dotvzeta1}.

We are now ready to apply Theorem \ref{thm:main:new:mat}. To that end, define
\begin{equation}
X := 
\begin{bmatrix}
x_1 \\
\vdots\\
x_n\\
\zeta_1\\
\vdots\\
\zeta_{n-1}
\end{bmatrix}\,,\quad
\Phi(t,X) := 
\begin{bmatrix}
\phi_1\\
\vdots\\
\phi_{n-1}
\end{bmatrix}\,.
\end{equation}
Further, define $\theta_i(\cdot):=\rho_{W_i}(\cdot)\rho_{B_i}(\cdot)$ from \rref{Bbnd}, \rref{rhoW} and notice that in particular, $\theta_i(0)=0$. Then, we have from \rref{eq:Lyap1}, \rref{dotvis2}, \rref{ch:dotvzeta1} and \rref{dotvksfinal}  that for all $i\in \{1,\ldots, n-1\}$ and almost all $(t,X)\in \mR\times {\cal B}(\Delta)^2$
\begin{eqnarray*}
\dot V_1 & \leq & -\rho\circ\kappa_n(\norm{X_n})\\
\dot V_2 & \leq & -c_n\norm{\Phi_{n-1}}^2 + \nu[ \,\theta_n(\norm{X_n}) + \sigma\circ\theta_n(\norm{X_n})\,] \\
\dot V_3 & \leq & -c_{n-1}\norm{\Phi_{n-2}}^2 + \nu[\,\norm{\Phi_{n-1}} + \theta_n(\norm{X_n})\, ]\\
\vdots\\
\dot V_{i+3} & \leq & -c_{n-i-1} \norm{\Phi_{n-i-2}(t,X)}^2 + \nu\left(\norm{\Phi_{n-i-1}(t,X)} + \norm{\Phi_{n-i}(t,X)} \right) \\
\vdots\\
\dot V_n  & \leq & -c_{2}\norm{\Phi_{1}(t,X)}^2 + \nu\left(\norm{\Phi_{2}(t,X)} + \norm{\Phi_{3}(t,X)} \right) \\
\vdots\\
\dot V_{n+i}  & \leq & -\norm{X_{n+i}}^2+ \nu\left(\theta_n(\norm{X_n}) + \sigma\circ\theta_n(\norm{X_n}) + \sum_{\ell=1}^{n-1} \norm{\Phi_{\ell}(t,X)}^{1/n-\ell}\,  \right)\\
\vdots\\
\dot V_{2n+i-1}  & \leq & -\gamma_{n-i}(\norm{X_{n-i}}) \norm{h_{n-i}(X_{n-i})}^2 +
 \nonumber \\ && 
\nu \left(\, \theta_n(\norm{X_n})  + \sigma\circ\theta_n(\norm{X_n}) + \norm{\Phi_{n-i}(t,X)}^{2/i} + \sum_{\ell=1}^{n-1} \norm{\Phi_{\ell}(t,X)}^{1/n-\ell} + \norm{X_{n+\ell}} \,  \right)\\[-0.8cm]
\vdots\\[0.5cm]
\dot V_{3n-2}  & \leq & -\gamma_1(\norm{X_{1}}) \norm{h_{1}(X_{1})}^2 +
 \nonumber \\ && 
\nu \left(\, \theta_n(\norm{X_n})  + \sigma\circ\theta_n(\norm{X_n}) + \norm{\Phi_{1}(t,X)}^{2} + \sum_{\ell=1}^{n-1} \norm{\Phi_{\ell}(t,X)}^{1/n-\ell} + \norm{X_{n+\ell}} \,  \right)
\end{eqnarray*}
Letting each of the bounds above be denoted by $Y_k(X,\Phi(t,X))$ with $k\in [1,\ldots, 3n-2]$ we see that under Assumption \ref{assongiab} each of these functions is bounded on compact sets of the space of $(X,\Psi = \Phi(t,X))$ uniformly in $t$, i.e., Assumption \ref{ass:mat:2} holds. Also, we see that each $Y_k(X,\Phi(t,X))$  is non-positive on the sets where all the previous bounds are identically zero (i.e., Assumption \ref{ass:mat:3} holds). Furthermore,  in view of \rref{ykappa}, in particular since $\kappa\in {\cal P}{\cal D}$, and since  $\gamma_i\in {\cal K}$ we see that the set where $Y_k(X,\Psi)\equiv 0$ is the origin, $\{X=0\}$ (i.e., Assumption \ref{ass:mat:4}). UGAS follows.
\qed 

\section{Proof of Theorem \ref{thm:main:new:mat}}
\label{sec:proof}
To prove the theorem we first need to establish the following claims. 
\begin{claim}\label{thm:mat:claim1}
Given $\delta>0$, there exists $\ep>0$ such that: 
\begin{description}
\item{(A):}
$ \quad
\left\{ \ (z,\psi) \in \HdD \times \Bmu \quad  \& \quad
Y_i(z,\psi) = 0   \right. $
$ \left. \forall i\in \{ 1,\ldots, j-1 \} \ \right\}$
\end{description}
implies
\begin{description}
\item{(B):}
$ \quad
\left\{ \ Y_j(z,\psi) \leq -\ep \ \right\} \ .
$
\end{description}
\end{claim}
\begin{claimproof}
We prove the claim by contradiction. Suppose that for each integer
$n$, there exist $(z_n,\psi_n)\in  \HdD\times \Bmu $ such that $
Y_i(z_n,\psi_n) = 0$ for all $i\in \{ 1,\ldots, j-1\}$, and
$Y_j(z_n,\psi_n) > - \frac{1}{n}$. By compactness of $\HdD\times \Bmu
$, the
continuity of $Y_j(\cdot,\cdot)$, and Assumption
\ref{ass:mat:3},
the sequence $(z_{n},\psi_{n})$ has an accumulation point
$(z_*,\psi_*)\in  \HdD\times \Bmu $ such that $Y_i(z_*,\psi_*)
= 0$ 
for all $i \in \{1,\ldots,j\}$.  By Assumption \ref{ass:mat:4},
this implies that $z_*=0$ which contradicts the fact that $z_* \in \HdD$.
\end{claimproof}

\begin{claim}\label{thm:mat:claim2}
Let $\ell\in \{2,\ldots, j \}$, $\widetilde\ep > 0$ and  a continuous
function $\widetilde Y_\ell:\mR^n \times \mR^{m} \rightarrow \mR$ be given.
Then the property
\begin{description}
\item{(1):} (A) implies (B) where 
\begin{description}
\item{(A):}
 $\quad \left\{ \ (z,\psi) \in \HdD \times \Bmu \quad \& \quad 
Y_i(z,\psi) = 0    \right. $
$ \left. \forall i\in \{ 1,\ldots, \ell-1\} \ \right\}$
\end{description}
\begin{description}
\item{(B):}
$\quad \{ \ \widetilde Y_\ell(z,\psi) \leq -\tilde \ep \ \} $
\end{description}
\end{description}
implies the property 
\begin{description}
\item{(2):}
there exists $K_{\ell-1} > 0 $ such that
\begin{description}
\item{(A):}
 $\quad \{ \ (z,\psi) \in \HdD \times \Bmu \quad \& \quad
Y_i(z,\psi) = 0 \quad \forall i\in \{ 1,\ldots, \ell-2\} \ \}$
\end{description}
 implies that 
\begin{description}
\item{(B):}
$\quad \{ \ K_{\ell-1} Y_{\ell-1}(z,\psi) + \widetilde Y_\ell(z,\psi) \leq -
\displaystyle \frac{\widetilde \ep}{2} \ \}$.
\end{description}
\end{description}
\end{claim}
\begin{claimproof}
By Assumption \ref{ass:mat:3} and ${\cal H}(\delta,\Delta) \subset {\cal B}(\Delta)$,  Property 2A implies that
$Y_{\ell-1}(z,\psi) \leq 0$.  Therefore Property 2A implies
\begin{equation}
 K_{\ell-1} Y_{\ell-1}(z,\psi)
+ \widetilde{Y}_{\ell}(z,\psi)
\leq \widetilde{Y}_{\ell}(z,\psi) \quad \forall K_{\ell-1} \geq 0 \ .
\end{equation}
Now if
$Y_{\ell-1}(z,\psi) =  0$ then, due to Property 1,
Property 2B holds for all $K_{\ell-1} \geq 0$ whenever Property 2A holds. 
We claim further that there exists $m>0$ such that
Property 2B holds whenever Property 2A holds and
$Y_{\ell-1}(z,\psi) > - m$.
Suppose not, i.e., for each integer $n$ there exists 
$(z_n,\psi_n)\in\HdD\times\Bmu$ such that 
$Y_{\ell-1}(z_n,\psi_n) > - \frac{1}{n}$ and 
\begin{equation}
\widetilde Y_{\ell}(z_{n},\psi_{n}) > \displaystyle -\frac{\widetilde
  \ep}{2} \ .
\label{eq:6AT}
\end{equation} 
Then, by compactness of $\HdD \times \Bmu$, continuity of 
${Y}_{\ell-1}$, and Assumption \ref{ass:mat:3},
the sequence $(z_{n},\psi_{n})$ has an accumulation point 
$(z_*,\psi_*)$ such that  
$Y_{\ell-1}(z_*,\psi_*) = 0$. 
But then from Property 1 we have that
$\widetilde Y_{\ell}(z^{*},\psi^{*})\leq -\tilde\ep$. 
By continuity of $\widetilde Y_{\ell}(\cdot,\cdot)$ 
this contradicts \rref{eq:6AT} when $n$ is large and associated
with a subsequence converging to the accumulation point.

It now follows from the continuity of $\widetilde{Y}_{\ell}$ and compactness of ${\cal H}(\delta,\Delta) \times {\cal B}(\mu)$ that we
can pick $K_{\ell-1}>0$ large enough to satisfy
\[
-m K_{\ell-1} + \max_{(z,\psi)\in\HdD\times\Bmu} \tilde Y_\ell(z,\psi)
 \, \leq \, - \frac{\widetilde\ep}{2}
\] 
then Property 2A implies Property 2B.
\end{claimproof}

We now use these two claims to prove the theorem. According
to Claim \ref{thm:mat:claim1}, Property 1 of Claim
\ref{thm:mat:claim2}
holds when $\ell=j$, $\widetilde{\epsilon} = \epsilon$ and
$\widetilde{Y}_{\ell}=Y_{j}$.  An application of Claim
\ref{thm:mat:claim2}
with these choices provides a value $K_{j-1}$ such that
Property 1 of Claim \ref{thm:mat:claim2} holds when
$\ell=j-1$, $\widetilde{\epsilon} = \epsilon/2$ and
$\widetilde{Y}_{\ell}=K_{j-1} Y_{j-1} + Y_{j}$.  Continuing with
this iteration, it follows that for each $\delta>0$ there
exists $\varepsilon>0$ and positive real numbers
$K_{i}$, $i=1,\ldots,j-1$ such that, for all $(z,\psi) \in \HdD \times \Bmu$
\begin{equation}\label{nine}
Z(z,\psi) := \sum_{i=1}^{j-1} K_{i} Y_{i}(z,\psi) + Y_{j}(z,\psi)
\leq - \frac{\ep}{2^{j-1}} \ .
\end{equation}

Next define the locally Lipschitz 
function $W:\mR \times \mR^{n} \rightarrow \mR$ as
\begin{equation}
W(t,x):= \sum_{i=1}^{j-1} K_{i} V_{i}(t,x) + V_{j}(t,x) \ .
\end{equation}
According to the conditions of the theorem and the discussion
above, we have that,
for almost all $(t,x) \in \mR \times \BD$, 
\begin{equation}
   |W(t,x)| \leq \mu \left( 1+ \sum_{i=1}^{j-1} K_{j-1} \right) =:
    \eta \ ,
\label{eq:eta_def}
\end{equation}
\begin{equation}
\dot W(t,x) \leq Z(x,\phi(t,x)) \ ,
\end{equation}
and, using \rref{bndonViPhii} together with \rref{nine} we obtain that for all $(t,x) \in \mR \times \HdD$,
\begin{equation}
  Z(x,\phi(t,x)) \leq - \frac{\varepsilon}{2^{j-1}} \ .
\end{equation}

Using Assumption \ref{ass:mat:1}, for each $r>0$ and $\rho>0$ there
exists $\Delta$ and $\delta$ such that
\begin{equation}
  |x_{\circ}| \leq r \quad \Longrightarrow \quad
  |x(t,t_{\circ},x_{\circ})| \leq \Delta \quad
  \forall t \geq t_{\circ}
\end{equation}
and
\begin{equation}
  |x(t_{1},t_{\circ},x_{\circ})| \leq \delta
\quad \Longrightarrow \quad |x(t,t_{\circ},x_{\circ})|
   \leq \rho \quad \forall t \geq t_{1} \ .
\end{equation}
Let $r$ and $\rho$ generate $\Delta$ and $\delta$ and
then let $\Delta$ and $\delta$ generate $\varepsilon$
and $\eta$ through the claims and definitions above.  Let 
\begin{equation}
  T > \frac{2^{j}\eta}{\varepsilon} \ .
\label{eq:Tchoice}
\end{equation}
We claim that 
\begin{equation}
  |x_{\circ}| \leq r \ , \ t \geq t_{\circ} + T
 \qquad \Longrightarrow \qquad
  |x(t,t_{\circ},x_{\circ})| \leq \rho \ .
\end{equation} 
Suppose not. It follows from the discussion above that
$x(t,t_{\circ},x_{\circ}) \in \HdD$ for all $t \in
[t_{\circ},t_{\circ}+T]$.
It then follows that, for almost all $t \in [t_{\circ},t_{\circ}+T]$,
\begin{equation}
  \dot{W}(t,x(t,t_{\circ},x_{\circ}))
   \leq - \frac{\varepsilon}{2^{j-1}} \ .
\end{equation}
Integrating and using (\ref{eq:eta_def}) we have
\begin{equation}
  T \frac{\varepsilon}{2^{j-1}} \leq 2 \eta
\end{equation}
which contradicts the choice of $T$ in (\ref{eq:Tchoice}).
\qed

\section{Conclusions}
\label{sec:concl}
In this paper we have presented a new tool for establishing
uniform attractivity of the origin when the origin is uniformly stable.
The tool involves the use of an arbitrary finite number of auxiliary
functions whose derivatives are simultaneously zero only at the origin.
This result generalizes Matrosov's theorem on uniform asymptotic stability of the origin for nonlinear time-varying systems.

As in other auxiliary-functions based methods {\em a la Lyapunov} it is not possible to give a general criterion on how to chose the auxiliary functions. However, we have given some intuition on how to use  our theorem by employing it to establish uniform global asymptotic stability of the origin for a controlled nonholonomic system. 


\bibliographystyle{ieeetr}
\bibliography{refs}


\appendix

\numberwithin{equation}{section}

\section{Useful identities involving the gains $g_i$}

We single out some important and not obvious identitites involving the gains $g_i$ defined in \rref{eq:fform} and which are used in the proof of Theorem \ref{thm:channels}.

\begin{itemize}
\item By simply spelling out the definition of $g_i$ we have that 
\begin{equation}
\label{2001}
g_i = \widetilde g_i \, g_{i+1} \quad \forall\, i\in [1,\ldots,n-2] \qquad\quad g_{n-1} = \tilde g_{n-1}\qquad\quad g_n > 0 \,.
\end{equation}
Furtheremore, using these identities one can show that for all  $i \in [1,\ldots,n-2]$,
\begin{subequations}\label{2000}
\begin{eqnarray}
\label{2000a}
g_i & = & [g_{n-1}\cdots g_i]^{\ \! 1/n-i} \ \widetilde g_{n-2}^{\ \! 1/n-i} \,\cdots\,\widetilde g_{i+j}^{\ \! n-i-j-1/n-i}\,\cdots\, \widetilde g_i^{\ \! n-i-1/n-i} \\
\label{2000b} 
g_{n-1} & = & \widetilde g_{n-1}\,.
\end{eqnarray}
\end{subequations}
For illustration, we have the following for $i=3$ and $n=6$.
\begin{eqnarray*}
g_3 &=& g_3^{1/3} \, g_3^{1/3} \, g_3^{1/3}\\
    &=& g_3^{1/3} \, (\widetilde g_3^{1/3} g_4^{1/3}) \, (\widetilde g_3^{1/3}\, g_4^{1/3} ) \\
&=& g_3^{1/3} \, (\widetilde g_3^{1/3} g_4^{1/3}) \, (\widetilde g_3^{1/3}\, [\widetilde g_4^{1/3}\, g_5^{1/3} ]\, )\\
  &=& (g_3^{1/3}\,g_4^{1/3}\,g_5^{1/3})\cdot \widetilde g_3^{1/3} \cdot \widetilde g_3^{1/3}\, \widetilde g_4^{1/3} \,.
\end{eqnarray*}
\item For the $n-1$ $\phi$'s  defined in \rref{ch:phis} which for convenience we index here as $\phi_i$ with $i\in [1,\ldots,n-2]$, and using the previous identities we have that 
\begin{subequations}\label{2002}
\begin{eqnarray}
\norm{g_i\, y_i} & = & \underbrace{\norm{[g_{n-1}\cdots g_i]}^{\ \! 1/n-i}\, \norm{y_i}^{1/n-i}}_{\dty \phi_i^{1/n-i}} \ \norm{y_i}^{-1+n-i/n-i} \ \norm{\widetilde g_{n-2}^{\ \! 1/n-i} \,\cdots\, \widetilde g_i^{\ \! n-i-1/n-i}}\nonumber\\[-6mm] \\[1mm]
\norm{g_{n-1}}\norm{ y_{n-1}} & = & \phi_{n-1}\,.
\end{eqnarray}
\end{subequations}
\end{itemize}

\end{document}